\documentstyle[amscd,12pt,leqno]{amsart}
\def\newsection{\section}
\numberwithin{equation}{section}
\renewcommand{\theequation}{\thesection.\arabic{equation}}
\newcommand{\g}{\frak}
\newcommand{\gtg}{{\g g}}
\def\Phim#1{\mathrel{\mathop{\kern0pt \Phi}\limits^#1}}
\def\Psim#1{\mathrel{\mathop{\kern0pt \Psi^*}\limits^#1}}
\def\on{\operatorname}
\newcommand{\eq}{\begin{eqnarray}}
\newcommand{\beq}{\begin{eqnarray}}
\newcommand{\beqn}{\begin{eqnarray*}}
\def\endeq{\end{eqnarray}}
\newcommand{\eqn}{\begin{eqnarray*}}
\def\endeqn{\end{eqnarray*}}
\newcommand{\nn}{\nonumber}
\def\C{{\Bbb C}}
\def\Z{{\Bbb Z}}
\def\Q{{\Bbb Q}}

\def\mod{{\operatorname {mod}}\,}
\def\id{{\operatorname {id}}}
\def\tr{{\operatorname {tr}}}
\newcommand{\univ}{{\on{univ}}}
\newcommand{\Ru}{R^\univ}

\def\End{{\operatorname {End}\,}}

\def\Pr{{\on P}}
\newcommand{\eps}{\varepsilon}

\def\hb{\hfill\break}

\newcommand{\gsl}{{\g{sl}}}

\newcommand{\uasl}{U'_q({\widehat{\g{sl}}_n})}

\newcommand{\bk}{{\bar k}}
\newcommand{\bA}{{\bar A}}
\newcommand{\U}{U_q({\g g})}
\renewcommand{\Im}{{\on{Im}}}

\theoremstyle{plain}
 \newtheorem{thm}{Theorem}[section]
 \newtheorem{prop}[thm]{Proposition}
 \newtheorem{lemma}[thm]{Lemma}
 \newtheorem{cor}[thm]{Corollary}
\newtheorem{claim}{Claim}

\theoremstyle{definition}
 \newtheorem{defn}[thm]{Definition}
 \newtheorem{conj}{Conjecture}

\theoremstyle{remark}
\newtheorem{rem}{Remark}

\newcommand{\Lemma}{\begin{lemma}}
\newcommand{\Def}{\begin{defn}}
\newcommand{\enDef}{\end{defn}}
\def\lan{\langle}
\def\ran{\rangle}
\def\Gt{{\g{t}}}
\def\Gg{{\g{g}}}
\def\lam{\lambda}
\def\Lam{\Lambda}
\def\al{\alpha}
\def\cl{{{\operatorname {cl}}}}
\def\clwt{\Gt^{*0}_\cl}
\def\Us{U'_q(\Gg)}
\def\Aut{{\operatorname{Aut}}}
\def\te{\tilde e}
\def\tf{\tilde f}
\newcommand{\proof}{\begin{pf}}
\newcommand{\QED}{\end{pf}}
\def\enlemma{\end{lemma}}

\newcommand{\wt}{{\on{wt}}}
\renewcommand{\max}{{\on {max}}}
\renewcommand{\min}{{\on {min}}}

\def\Definition{\begin{defn}}
\def\Thm{\begin{thm}}
\def\enthm{\end{thm}}
\def\Corollary{\begin{cor}}
\def\encorollary{\end{cor}}
\def\Cor{\begin{cor}}
\def\encor{\end{cor}}
\def\Prop{\begin{prop}}
\def\enprop{\end{prop}}

\def\owari{\end{document}}
\newcommand{\ol}{\overline}
\newcommand{\Rn}{R^{{\on{nor}}}}
\newcommand{\Tr}{\cal{T}\!r}
\newcommand{\Hom}{\on{Hom}}
\newcommand{\p}{{p^*}}
\newcommand{\To}{\longrightarrow}
\newcommand{\ip}[1]{\{#1\}}
\newcommand{\m}{{\g{m}}}
\newcommand{\Prod}{\mathop{\displaystyle\prod}\limits}
\newcommand{\ake}{\hspace{5pt}}
\def\Gb{{\g{b}}}
\newcommand{\Usb}{U'_q(\Gb)}

\title[Finite-dimensional representations
of quantum affine algebras]
{Finite-dimensional representations
of quantum affine algebras}

\author{Tatsuya Akasaka and Masaki Kashiwara}
\address[]{Department of Mathematics and
Research Institute for Mathematical Sciences
Kyoto University, Kyoto 606, Japan}

\begin{document}
\maketitle

\begin{abstract}
We present a conjecture on the irreducibility of the
tensor products of fundamental representations
of quantized affine algebras. 
This conjecture implies in particular that the irreducibility of
the tensor products of fundamental representations
is completely described by the poles of $R$-matrices.
The conjecture is proved in the cases of type $A^{(1)}_n$ and $C^{(1)}_n$.
\end{abstract}




\setcounter{section}{-1}

\newsection{Introduction}

In this paper we study finite-dimensional representations
of quantum affine algebras. It is known that any
finite-dimensional irreducible representation is isomorphic to the
irreducible subquotient of a tensor product
$\otimes_\nu V(\varpi_{i_\nu})_{a_\nu}$ containing
the highest weight (Drinfeld \cite{Drin}, Chari-Pressley \cite{CP}).
Here $V(\varpi_i)$ is the fundamental representation
corresponding to the fundamental weight $\varpi_i$ and $a_\nu$ are
spectral parameters.
Moreover $\{ (\varpi_{i_\nu};a_\nu)\}_\nu$ is uniquely determined
up to permutation. This gives a parameterization of the isomorphic
classes of finite-dimensional irreducible representations.

However it is not known for example 
what is the character of those irreducible representations 
except the complete result for $A^{(1)}_1$ (\cite{CP}) and 
some other results due to Chari-Pressley (\cite{CP,CP1,CP2}). 
We have even not known 
when $\otimes V(\varpi_{i_\nu})_{a_\nu}$ itself is irreducible.

In this paper we propose a conjecture on the irreducibility of 
$\otimes_\nu V(\varpi_{i_\nu})_{a_\nu}$
and prove this conjecture for $A^{(1)}_n$ and $C^{(1)}_n$. 

For $x,y\in\C(q)$, let us denote $x\le y$ if $x/y$ does not
have a pole at $q=0$.
We denote by $u_{i}$ the highest weight
vector of $V(\varpi_{i})$.

\begin{conj}\label{conj:main}  \hb
(1)\quad If $a_{1} \leq \cdots \leq a_{N}$, then 
$V(\varpi_{i_{1}})_{a_{1}}\otimes\cdots\otimes 
V(\varpi_{i_{N}})_{a_{N}}$ 
is generated by $u_{i_{1}}\otimes\cdots\otimes u_{i_{N}}$ 
as a $\Us$-module.
\hb
(2)\quad If $a_{1} \geq \cdots \geq a_{N}$,
then any non-zero $\Us$-submodule of 
$V(\varpi_{i_{1}})_{a_{1}}\otimes\cdots\otimes V(\varpi_{i_{N}})_{a_{N}}$ 
contains $u_{i_{1}}\otimes\cdots\otimes u_{i_{N}}$.
\end{conj}
\newcounter{old}
\setcounter{old}{\value{conj}}

Here $\Us$ is the quantum affine algebra without 
derivation (see \S\ref{sec:aff}).

This conjecture implies in particular the following 
consequences.

\begin{claim}
If $a_{1} \leq a_{2}$,
then the normalized $R$-matrix
$$\Rn_{i,j}(x,y)
: V(\varpi_{i})_x\otimes V(\varpi_{j})_y\to V(\varpi_{j})_y
\otimes V(\varpi_{i})_x
$$ does not have a pole
at $(x,y)=(a_1,a_2)$.
\end{claim}
Here $\Rn_{i,j}(x,y)$ is so normalized that it sends 
$u_{i}\otimes u_{j}$ to
$u_{j}\otimes u_{i}$.

\begin{claim}
$V(\varpi_{i_1})_{a_1}\otimes\cdots\otimes V(\varpi_{i_N})_{a_N}$
is irreducible if and
only if the $R$-matrix
$$
\Rn_{i_\nu,i_\mu}(x,y)
: V(\varpi_{i_\nu})_x\otimes V(\varpi_{i_\mu})_y\to V(\varpi_{i_\mu})_y
\otimes V(\varpi_{i_\nu})_x
$$
does not have a pole at $(x,y)=(a_\nu,a_\mu)$ for any 
$1\le \nu,\,\mu\le N\,(\nu\not =\mu)$.
\end{claim}

\begin{claim}
Assume that $\Rn_{i_\nu,i_\mu}(x,y)$ has no pole at
$(x,y)=(a_\nu,a_\mu)$ for any
$1\le \mu<\nu\le N$.
Then the submodule generated by
$u_{i_1}\otimes\cdots\otimes u_{i_N}$
is an irreducible submodule of
$V(\varpi_{i_1})_{a_\nu}\otimes\cdots\otimes V(\varpi_{i_N})_{a_N}$.
Conversely, any finite-dimensional irreducible integrable module
is obtained in this way.
\end{claim}

\begin{claim}
If $M$ and $M'$ are irreducible finite-dimensional integrable
$\Us$-modules,
then
$M\otimes M'_z$ is an irreducible $\Us$-module 
except for finitely many $z$.
\end{claim}

The plan of the paper is as follows.
In \S \ref{Notations},
we fix notations and explain the results
used later.
We announce non published results but they can be directly
checked for the $A^{(1)}_n$ and $C^{(1)}_n$ cases.
In \S \ref{Conjecture}, we announce the main conjecture
and discuss its consequences.
In \S \ref{Aux.Conj}, we reduce the main conjecture
to another auxiliary conjecture,
which will be proved in the case $A^{(1)}_n$ and 
$C^{(1)}_n$ in \S \ref{AC}.
In the appendix, we shall calculate the explicit
form of the normalized $R$-matrices and the universal $R$-matrices between
fundamental representations of $A^{(1)}_n$ and $C^{(1)}_n$.

The authors are grateful to K. Takemura 
for his helpful comments on this work.

\newsection{Notations}\label{Notations}
\subsection{Quantized affine algebras}\label{sec:aff}

 Let $(a_{ij})_{i,j\in I}$ be a generalized Cartan matrix of affine type.
We choose a $\Q$-vector space $\Gt$ of dimension $\sharp I+1$
and simple roots $\alpha_i\in\Gt^*$ and simple coroots $h_i\in\Gt$
such that
$\lan h_i,\alpha_j\ran=a_{ij}$.
We assume further that $\alpha_i$ and $h_i$ are linearly independent.
Set $Q=\sum_i\Z\alpha_i$ and $Q^\vee=\sum_i\Z h_i$.
Let $\delta=\sum a_i\alpha_i$ be the smallest positive
imaginary root and let
$c=\sum a_i^\vee h_i\in Q^\vee$ be the center.
Set $\Gt^*_\cl=\Gt^*/\Q\delta$ and
let $\cl:\Gt^*\to \Gt^*_\cl$ be the projection.
We set $\Gt^{*}{}^{0}=\{\lam\in\Gt^*;\lan c,\lam\ran=0\}$
and $\clwt=\cl(\Gt^{*}{}^{0})$.

We take a non-degenerate symmetric bilinear form $(\cdot,\cdot)$
on $\Gt^*$
such that
\eqn
&&\lan h_i,\lam\ran={2(\alpha_i,\lam)\over(\alpha_i,\alpha_i)}
\quad\mbox{for any $i\in I$ and $\lam\in\Gt^*$.}
\endeqn
We normalize it by
\eq
&&\lan c,\lam\ran=(\delta,\lam)\quad\text{for any $\lam\in \Gt^*$.}
\endeq
We identify sometimes $\Gt$ and $\Gt^*$ by this symmetric form.

Let us take a (weight) lattice $P\subset\Gt^*$
such that
$\alpha_i\in P$ and $h_i\in P^*$ for every $i\in I$.
We assume further that $P$ contains $\Lam_i$ satisfying
$\lan h_j,\Lam_i\ran=\delta_{ij}$
and that $P\cap \Q\delta=\Z\delta$.
We set $P_\cl=P/\Z\delta\subset \Gt^*_\cl$,
$P^0=\{\lam\in P\,;\,\lan c,\lam\ran=0\}\subset\Gt^{*}{}^{0}$,
and $P^0_\cl=\cl(P^0)\subset \clwt$.
Note that the dual lattice of $Q^\vee$ coincides
with $P_\cl\cong\oplus_{i\in I}\Z\cl(\Lam_i)$.

Let $\gamma$ be the smallest positive integer such that
\eq\label{def:gamma}
&&\gamma(\alpha_i,\alpha_i)/2\in\Z\quad\text{for any $i\in I$.}
\endeq
Then the quantized affine algebra
$\U$ is the algebra over $k=\Q(q^{1/\gamma})$ generated by the symbols
$e_i,\,f_i$$(i\in I)$ and $q(h)\,(h\in\gamma^{-1} P^*)$
satisfying the following defining relations.
\begin{enumerate}
\item $q(h)=1$ for $h=0$.

\item $q(h_1)q(h_2)=q(h_1+h_2)$ for $h_1,h_2\in\gamma^{-1} P^*$.

\item For any $i\in I$ and $h\in\gamma^{-1} P^*$,
\begin{eqnarray*}
q(h)e_iq(h)^{-1}&=&q^{\langle h,\al_i\rangle}e_i\quad \mbox{and}\\
q(h)f_iq(h)^{-1}&=&q^{-\langle h,\al_i\rangle}f_i\,.
\end{eqnarray*}

\item\label{even} $\lbrack e_i,f_j\rbrack
=\delta_{ij}\dfrac{t_i-t_i^{-1}}{q_i-q_i^{-1}}$
for $i,j\in I$. Here $q_i=q^{(\al_i,\al_i)/2}$ and
$t_i=q(\frac{(\al_i,\al_i)}{2}h_i)$.

\item (Serre relations) For $i\not= j$,
\begin{eqnarray*}
&&\sum^b_{k=0}(-1)^ke^{(k)}_ie_je^{(b-k)}_i=\sum^b_{k=0}(-1)^kf^{(k)}_i
f_jf_i^{(b-k)}=0.
\end{eqnarray*}
Here $b=1-\langle h_i,\al_j\rangle$ and
$$
\begin{array}{ll}
e^{(k)}_i=e^k_i/\lbrack k\rbrack_i!\ ,& f^{(k)}_i=f^k_i/\lbrack k\rbrack_i!\ ,\\
\lbrack k\rbrack_i=(q^k_i-q^{-k}_i)/(q_i-q^{-1}_i)\ ,
&\lbrack k\rbrack_i!=\lbrack 1\rbrack_i\cdots \lbrack k\rbrack_i\,.
\end{array}
$$
\end{enumerate}

We denote by $\Us$ the subalgebra of $\U$ generated by
$e_i,\,f_i\,(i\in I)$ and $q(h)\,(h\in\gamma^{-1} Q^\vee)$.

In this paper we consider only $\Us$.
A $\Us$-module $M$ is called {\sl integrable}
if $M$ has the weight decomposition
$M=\oplus_{\lam\in P_\cl}M_\lam$ where
$M_\lam=\{u\in M;q(h)u=q^{\lan h,\lam\ran}u\}$,
and if $M$ is $\U_i$-locally finite
(i.e. $\dim\U_iu<\infty$ for every $u\in M$)
for every $i\in I$.
Here $\U_i$ is the subalgebra generated
by $e_i$, $f_i$ and $t_i$.

We use the coproduct $\Delta$ of $\U$
given by
\eq
\Delta(q(h))&=&q(h)\otimes q(h)\,,\\
\Delta(e_i)&=&e_i\otimes t_i^{-1}+1\otimes e_i\,,\\
\Delta(f_i)&=&f_i\otimes 1+t_i\otimes f_i\,,
\endeq
so that the lower crystal bases behave well
under the corresponding tensor products
(\cite{Cry}).

\subsection{Finite-dimensional representations}
Let $W\subset\Aut(\Gt^*)$ be the Weyl group,
and let $l:W\to \Z$ be the length function.
Since $\delta$ is invariant by $W$,
we have the group homomorphism
$\cl_0:W\to\Aut(\clwt)$.
Let $W_\cl\subset \Aut(\clwt)$ be the image of $W$
by $\cl_0$.
Then $W_\cl$ is a finite group.
Let us take $i_0\in I$ such that $W_\cl$ is generated by $\cl_0(s_i)$
($i\in I_0=I\setminus\{i_0\}$) and that
$a_{i_0}^\vee=1$. 
Such an $i_0$ is unique up to Dynkin diagram automorphism.
Hereafter we write $0$ instead of $i_0$. We have $(\alpha_0,\alpha_0)=2$.

Let us denote by $W_0$ the subgroup of $W$ generated by
$s_i$ ($i\in I_0=I\setminus\{0\}$).
Then $W_0$ is isomorphic to $W_\cl$.
The kernel of $W\to W_\cl$
is the commutative group
$\{t(\xi);\xi\in Q_\cl\cap Q^\vee_\cl\}$.
Here
$Q_\cl=\cl(Q)=\sum_{i\in I}\Z\cl(\alpha_i)$ and $Q^\vee_\cl=\cl(Q^\vee)
=\sum_{i\in I_0}\Z\cl(h_i)$
and $t(\xi)$ is the automorphism of $\Gt^*$ given by
\eqn
&&t(\xi)(\lam)
=\lam+(\delta,\lam)\xi'-(\xi',\lam)\delta
-{(\xi',\xi')\over 2}(\delta,\lam)\delta
\endeqn
for $\xi'\in\Gt^*$ such that $\cl(\xi')=\xi$.

The following lemma is well-known.
\Lemma\label{lemma:21}
Let $\xi\in Q_\cl\cap Q^\vee_\cl$ and $w\in W_0$.
\hb
$\on{(i)}\quad$
If $\xi$
is dominant $($with respect to $I_0$$)$, then we have
$$l(w\circ t(\xi))=l(w)+l(t(\xi)).$$
$\on{(ii)}\quad$
If $\xi$
is regular and dominant, then we have
$$l(t(\xi)\circ w)=l(t(\xi))-l(w).$$
\enlemma

Let us choose $i_1$ such that
$W_\cl$ is generated by $\cl_0(s_i)$
($i\in I\setminus\{i_1\}$) and that
$a_{i_1}=1$.
For any $z\in k\setminus\{0\}$, let
$\psi(z)$ be the automorphism of $\Us$ given by
\eqn
\psi(z)(e_i)&=&z^{\delta_{i,i_1}}e_i\,,\\
\psi(z)(f_i)&=&z^{-\delta_{i,i_1}}f_i\,,\\
\psi(z)(q(h))&=&q(h)\,.
\endeqn
For a $\Us$-module $M$, let $M_z$ be the $\Us$-module
with $M$ as its underlying $k$-vector space
and with $\Us\overset{\psi(z)}\longrightarrow\Us\To\End(M)$
as the action of $\Us$.
Then $M\mapsto M_z$ is a functor satisfying
$(M\otimes N)_z\cong M_z\otimes N_z$.
This definition extends to the case $z\in K\setminus\{0\}$
for a field extension $K\supset k$.

If $M$ is a finite-dimensional integrable $\Us$-module, then
the weights of $M$ are contained in $P^0_\cl$.

\subsection{Fundamental representations}
We set $\varpi_i=\cl(\Lam_i-a_i^\vee\Lam_0)$ for
$i\in I_0$.
Then $(\varpi_i)_{i\in I_0}$ forms a basis of $P_\cl^0$.
We call $\varpi_i$ a {\em fundamental weight} (of level $0$).

For $i\in I_0$, there exists
an irreducible integrable $\Us$-module $V(\varpi_i)$
satisfying the following properties.

\begin{enumerate}
\item The weights of $V(\varpi_i)$
are contained in the convex hull of $W_\cl\varpi_i$.

\item $\dim V(\varpi_i)_{\varpi_i}=1$.

\item For any $\mu\in W_\cl\varpi_i\subset P_\cl^0$, 
we can associate a non-zero vector $u_\mu$
of weight $\mu$ such that
\eqn
u_{s_i\mu}&=&
\left\{\begin{array}{ll}
f_i^{(\lan h_i,\mu\ran)}u_\mu&\mbox{if $\lan h_i,\mu\ran\ge0$,}\\
e_i^{(-\lan h_i,\mu\ran)}u_\mu&\mbox{if $\lan h_i,\mu\ran\le0$.}
\end{array}
\right.
\endeqn
\end{enumerate}
Then $V(\varpi_i)$ is unique up to an isomorphism.
Moreover $V(\varpi_i)$ has a global crystal base.
We call $V(\varpi_i)$ a {\em fundamental representation}.
Then $V(\varpi_i)$ has a non-degenerate symmetric bilinear form
$(\,\cdot\,,\,\cdot\,)$ such that
${}^te_i=f_i$ and ${}^tq(h)=q(h)$.
Hence the duality is given as follows.
Let $w_0$ be the longest element of $W_0$.
Then for $i\in I_0$ there exists $i^*\in I_0$ such that
$$\varpi_{i^*}=-w_0\varpi_i.$$
(Remark that $i\mapsto i^*$ with $0^*=0$
gives a Dynkin diagram automorphism.)\hb
Then the right dual of $V(\varpi_i)$ is
$V(\varpi_{i^*})_p$ with the duality morphisms:
\eq\label{duality}
k\to
V(\varpi_{i^*})_\p\otimes V(\varpi_i)
\quad\text{and }
\quad
V(\varpi_i)\otimes V(\varpi_{i^*})_\p
\to k
\endeq
with $\p={(-1)^{\lan\rho^\vee,\delta\ran}q^{(\rho,\delta)}}$.
Here $\rho$ and $\rho^\vee$ are defined by:
$\lan h_i,\rho\ran=1$ and $\lan \rho^\vee,\alpha_i\ran=1$ for every $i\in I$.
Usually $(\rho,\delta)=\sum_{i\in I}a^\vee_i$ 
is called the dual Coxeter number and 
$\lan\rho^\vee,\delta\ran=\sum_{i\in I}a_i$ the Coxeter number.

Let $m_i$ be a positive integer such that
$$W(\Lam_i-a^\vee_i\Lam_0)=(\Lam_i-a^\vee_i\Lam_0)+\Z m_i\delta\,.$$
We have $m_i=(\alpha_i,\alpha_i)/2$
in the case where $\Gg$ is the dual of an untwisted affine algebra,
and $m_i=1$ in the other cases.

Then 
for $z,\,z'\in K^*$, we have 
\eq
&&V(\varpi_i)_z\cong V(\varpi_i)_{z'}
\quad
\mbox{if and only if $z^{m_i}=z'{}^{m_i}$.}
\endeq

Hence we set
$$V(\varpi_i;z^{m_i})=V(\varpi_i)_z\,.$$

The following theorem is announced by Drinfeld (\cite{Drin})
in Yangian case,
and its proof is given by Chari-Pressley (\cite{CP1,CP2}).


\Thm\label{drin}
Let $K\supset k$ be an algebraically closed field and
let $M$ be an irreducible finite-dimensional $\Us_K$-module.
Then there exist
$i_1,\ldots,i_N\in I_0$
and $z_1,\ldots z_N\in K\setminus\{0\}$
such that
$M$ is isomorphic to a unique irreducible subquotient of
$V(\varpi_{i_1};z_1)\otimes\cdots\otimes V(\varpi_{i_N};z_N)$
containing the weight $\sum_{\nu=1}^{N}\varpi_{i_\nu}$.
Moreover, $\{(i_1;z_1),\ldots,(i_N;z_N)\}$
is unique up to permutations.
\enthm

\Def\label{component}
We call $V(\varpi_{i_\nu};z_\nu)$ 
a component of $M$.
\enDef

\subsection{Extremal vectors}

We say that a crystal $B$ over $\Us$
is {\it a regular crystal} 
if, for any $J\underset{\not=}{\subset}I$,
$B$ is isomorphic to
the crystal associated with
an integrable $U_q(\Gg_J)$-module.
Here $U_q(\Gg_J)$ is the subalgebra of $\Us$ generated by
$e_i,f_i$ and $t_i$ ($i\in J$).
This condition is equivalent to saying that 
the same assertion holds for any $J\underset{\not=}{\subset}I$
with two elements (see \cite[Proposition 2.4.4]{(KMN)^2}).

By \cite{modified},
the Weyl group $W$ acts on any regular crystal.
This action $S$ is given by
\eqn
&&S_{s_i}b=
\begin{cases}
\tf_i^{\lan h_i,\wt(b)\ran}b
&\mbox{if $\lan h_i,\wt(b)\ran\ge 0$}\\
\te_i^{-\lan h_i,\wt(b)\ran}b
&\mbox{if $\lan h_i,\wt(b)\ran\le 0$.}
\end{cases}
\endeqn

A vector $b$ of a regular crystal $B$
is called $i$-extremal if
$\te_ib=0$ or $\tf_ib$=0.
We call $b$ an {\em extremal} vector
if $S_wb$ is $i$-extremal for any $w\in W$ and $i\in I$.

\Lemma\label{loop}
For any $\lam$, $\mu\in \clwt$ in the same $W_\cl$-orbit, 
we can find $i_1,\ldots,i_N$
$\in I$
such that
\beqn
&&\mu=s_{i_N}\cdots s_{i_1}\lam,\\
&&\lan h_{i_k},s_{i_{k-1}}\cdots s_{i_1}\lam\ran>0
\ \hbox{ for any $1\le k\le N$.}
\endeqn
\enlemma

\proof
It is enough to prove the statement above
for a regular integral anti-dominant 
(with respect to $I_0$) weight
$\lam$
and the dominant weight $\mu\in W\lam$.
We may assume further $\lam\in Q_\cl\cap Q_\cl^\vee$.
Let $w_0$ be the longest element of
$W_0$.
\hb
By Lemma \ref{lemma:21}, we have
\eqn
l(t(\lam))&=&l(t(-\lam))\\
&=&l(w_0)+l(t(-\lam)w_0)\\
&=&l(w_0)+l(w_0t(\lam)).
\endeqn
Take a reduced expression
$w_0t(\lam)=s_{i_N}\cdots s_{i_1}$.
Then for $1\le k\le N$ we have
$l(t(\lam)s_{i_1}\cdots s_{i_k})=
l(t(\lam))-k$ and hence
$t(\lam)s_{i_1}\cdots s_{i_{k-1}}\alpha_k$
is a negative root.
Since it is equal to
$s_{i_1}\cdots s_{i_{k-1}}\alpha_k
-(\lam,s_{i_1}\cdots s_{i_{k-1}}\alpha_k)\delta$
and $s_{i_1}\cdots s_{i_{k-1}}\alpha_k$ is a positive root,
we conclude
\beqn
(\lam,s_{i_1}\cdots s_{i_{k-1}}\alpha_k)>0.
\endeqn
On the other hand
we have the equality
$s_{i_N}\cdots s_{i_1}\lam
=w_0t(\lam)\lam=w_0\lam$
in $\clwt$.
Hence it is equal to $\mu$.
\QED

For a regular crystal $B$, $b\in B$ and $i\in I$,
let us denote by 
$\te_{i}^\max b$ the $i$-highest weight vector in the $i$-string 
containing $i$.
Namely we have
$$\te_{i}^\max b=\te_i^{\eps_i(b)}b.$$

\Lemma\label{ext}
Let $B$ be a finite regular crystal with level $0$
$($with weight in $P_\cl^0$$)$.
\begin{enumerate}
\item
For $b\in B$, there are $i_1,\cdots,i_N\in I$ such that
$\te_{i_N}^\max\cdots\te_{i_1}^\max b$ is an extremal vector.
\item
Any vector in the $W$-orbit of an extremal vector $b$ of $B$
is written in the form
$\te_{i_N}^\max\cdots\te_{i_1}^\max b$.
\end{enumerate}
\enlemma
\proof
Let us set
$F_l=
\{\te_{i_l}^\max\cdots\te_{i_1}^\max b\,;\,
i_1,\cdots,i_l\in I\}$,
$F=\bigcup_{l\ge0}F_l$.
Replacing $b$ with $b'\in F$
with maximal $(\wt(b'),\wt(b'))$,
we may assume from the beginning that
$(\wt(b'),\wt(b'))\le(\wt(b),\wt(b))$
for any $b'\in F$.
Since $(\wt(b'),\wt(b'))\ge(\wt(b),\wt(b))$,
we have $(\wt(b'),\wt(b'))=(\wt(b),\wt(b))$
for any $b'\in F$, and hence 
any $b'\in F$ is $i$-extremal for every $i\in I$.
Moreover the weight of $b'$ is in the $W_\cl$-orbit of $\wt(b)$.
Then for any weight $\mu$ of $F$ and
$i$ such that $\lan h_i,\mu\ran\le0$,
$S_{s_i}$ sends injectively $F_\mu$ to $F_{s_i\mu}$.
Hence $\sharp(F_\mu)\le\sharp(F_{s_i\mu})$, and
Lemma \ref{loop} asserts that they must be equal.
Therefore $S_i:F_\mu\to F_{s_i\mu}$
is bijective.
This shows that $F$ is stable by all $S_{s_i}$.
Thus we have (1) and (2).
\QED

\Lemma\label{exttensor:cry}
Let $B_1$ and $B_2$ be
two finite regular crystals.
Let $b_1$ and $b_2$ be vectors in $B_1$ and $B_2$,
respectively.
\begin{enumerate}
\item
If $b_1$ and $b_2$ are extremal vectors
and if their weights are in the same Weyl chamber,
then $b_1\otimes b_2$ is extremal.
\item
Conversely if $b_1\otimes b_2$ is extremal,
then
$b_1$ and $b_2$ are extremal vectors
and their weights are in the same Weyl chamber.
\end{enumerate}
\enlemma

\proof
(1) is obvious
because $S_w(b_1\otimes b_2)=S_wb_1\otimes S_wb_2$
under this condition.

We shall prove (2).
Since $\te_{i_1}^\max\cdots\te_{i_N}^\max(b_1\otimes b_2)
=\te_{i_1}^\max\cdots\te_{i_N}^\max b_1\otimes b'_2$
for some $b'_2\in B_2$,
the preceding lemma implies that $b_1$ is extremal.
Similarly $b_2$ is extremal.
It remains to prove that
$\wt(b_1)$ and $\wt(b_2)$
are in the same Weyl chamber.
Let us show first that $\wt(b_1\otimes b_2)$ and $\wt(b_1)$
are in the same Weyl chamber.
We may assume without loss of generality that 
$\wt(b_1\otimes b_2)$ is dominant (with respect to $I_0$).
Then $\te_i(b_1\otimes b_2)=0$ for every $i\in I_0$.
Hence $\te_ib_1=0$. Hence $\wt(b_1)$
is dominant.
Hence $\wt(b_1\otimes b_2)$ and $\wt(b_1)$
are in the same Weyl chamber.
Similarly $\wt(b_1\otimes b_2)$ and $\wt(b_2)$
are in the same Weyl chamber.
Thus $\wt(b_1)$ and $\wt(b_2)$
are in the same Weyl chamber.
\QED

\Definition\label{unext}
We say that a finite regular crystal $B$ is simple
if $B$ satisfies
\begin{enumerate}
\item There exists $\lam\in P_\cl^0$ such that
the weights of $B$ are in the convex hull of $W_\cl\lam$.
\item $\sharp(B_\lam)=1$.
\item The weight of any extremal vector is in $W_\cl\lam$.
\end{enumerate}
\enDef

\Prop
The crystal graph of the fundamental representations is simple.
\enprop

The proof will be given elsewhere.
However we can easily check this for the $A^{(1)}_n$ and $C^{(1)}_n$ cases.

\Lemma
A simple crystal $B$ is connected.
\enlemma

\proof
In fact, any vector is connected with an extremal vector
by Lemma \ref{ext}.
\QED

\Lemma
The tensor product of simple crystals is also simple.
\enlemma

\proof 
This immediately follows from Lemma \ref{exttensor:cry}.
\QED

\Prop
Let $M$ be a finite-dimensional integrable $\Us$-module with a crystal base
$(L,B)$.
Assume the following conditions.
\addtocounter{equation}{1}
\hb
$(\theequation)$\quad $B$ is connected.
\addtocounter{equation}{1}
\hb
$(\theequation)$\quad 
There exists a weight $\lam\in P^0_\cl$ such that
$\dim(M_\lam)=1$.\hb
Then $M$ is irreducible.
\enprop

\proof
We shall show first
that $M_\lam$ generates $M$.
Set $N=\Us M_\lam$ and
$\bar N=(L\cap N)/(qL\cap N)\subset L/qL$.
Then $\bar N$ is invariant by $\te_i$ and $\tf_i$.
Hence $\bar N$ contains $B$,
and Nakayama's lemma asserts that
$N=M$.
By duality, any non-zero submodule of $M$ contains $M_\lam$.
Therefore $M$ is irreducible.
\QED

\Cor
A finite-dimensional $\Us$-module with a simple crystal base
is irreducible.
\encor
\Cor\label{cor:irr}
For $i_1,\ldots, i_N\in I_0$,
$V(\varpi_{i_1})\otimes\cdots\otimes V(\varpi_{i_N})$
is irreducible.
\encor

We define similarly an extremal vector of
an integrable $\Us$-module.

\Def
Let $v$ be a weight vector of an integrable $\Us$-module.
We call $v$ {\em extremal}
if the weights of $\Us v$ are contained in the convex hull of
$W\wt(v)$.
\enDef

When the weight of $v$ is of level $0$ and dominant
(with respect to $I_0$),
$v$ is extremal if and only if
$\wt(\Us v)\subset\wt(v)+\sum_{i\in I_0}\Z_{\le0}\cl(\alpha_i)$.
In this case, we call $v$ a {\em dominant extremal vector}.

Since the following proposition is not used in this paper,
the proof will be given elsewhere.

\Prop
Let $v$ be a weight vector of an integrable $\Us$-module.
The following two conditions are equivalent.
\begin{enumerate}
\item
$v$ is an extremal vector.
\item
We can associate a vector $v_w$ of weight $w\wt(v)$
to each $w\in W$ satisfying the following properties:
\begin{enumerate}{\roman{enumii}}
\item
$v_w=v$ if $w=e$,
\item
If $i\in I$ and $w\in W$ satisfy
$\lan h_i,w\wt(v)\ran\ge0$, then
$e_iv_w=0$ and $v_{s_iw}=f_i^{(\lan h_i,w\wt(v)\ran)}v_w$,
\item
If $i\in I$ and $w\in W$ satisfy
$\lan h_i,w\wt(v)\ran\le0$, then
$f_iv_w=0$ and $v_{s_iw}=e_i^{-(\lan h_i,w\wt(v)\ran)}v_w$.
\end{enumerate}
\end{enumerate}
\enprop

The implication (1)$\Rightarrow$(2) is obvious.

%

\smallskip
Let us denote by
$\Usb$ the subalgebra of $\Us$ generated by
$t_i$ and $e_i$ ($i\in I$).

\Prop\label{prop:bs}
Let $M$ be a finite-dimensional integrable $\Us$-module.
Then any $\Usb$-submodule of $M$ is a
$\Us$-submodule.
\enprop
\proof
Let $N$ be a $\Usb$-submodule.
For any pair of weights $\lam$ and $\mu$ conjugate by $W_\cl$,
there exist $i_1,\ldots,i_l$ such that
$m_k=-\lan h_{i_{k}},s_{i_{k-1}}\cdots s_{i_1}\lam\ran>0$
and $\mu=s_{i_{l}}\cdots s_{i_1}\lam$ by Lemma \ref{loop}.
Then $e_{i_l}^{m_l}\cdots e_{i_1}^{m_1}$
sends injectively $N_\lam$ to $N_\mu$.
Hence we have $\dim N_\lam\le \dim N_\mu$.
Thus we obtain $\dim N_\lam=\dim N_\mu$.
Then the proposition follows from the following lemma.
\QED

\Lemma
Let $M$ be a finite-dimensional integrable $U_q(\gsl_2)$-module
and let $N$ be a vector subspace of $M$ stable by $e$ and $t$.
If $\dim N_\lam=\dim N_{s(\lam)}$ for any $\lam$
$($$s$ is the simple reflection$)$,
then $N$ is a $U_q(\gsl_2)$-submodule.
\enlemma
\proof
Any $u\in N_\lam$ can 
be written
\eqn
&&u=\sum_{n}f^{(n)}v_n
\endeqn
with $ev_n=0$.
Here $n$ ranges over $\{n\in \Z_{\ge0};n+\lan h,\lam\ran\ge 0\}$.

Let us prove $U_q(\gsl_2)v_n\subset N$
by the descending induction on $c=\lan h,\lam\ran$.
We have $eu=\sum_{n}[1+c+n]f^{(n-1)}v_n$.
Hence the induction hypothesis implies
$U_q(\gsl_2)v_n\subset N$ for $n>0$.
Hence we may assume that $eu=0$,
and then $c\ge 0$.
The surjectivity of $e^{c}:N_{s\lam}\to N_\lam$ implies the existence
of $w\in N_{s\lam}$ such that $u=e^{(c)}w$.
Then $fw=0$ and $U_q(\gsl_2)u=U_q(\gsl_2)w$ 
is generated by $\{e^n w;n\ge0\}\subset N$.
\QED


\Lemma\label{ext:mod}
Let $M_1$ and $M_2$ be finite-dimensional
$\Us$-modules and
let $v_1$ and $v_2$ be
non-zero weight vectors
of $M_1$ and $M_2$.
If $v_1\otimes v_2$ is extremal,
then
$v_1$ and $v_2$ are extremal and
their weights are in the same Weyl chamber
$($in $\clwt$$)$.
\enlemma
\proof
We may assume that $\wt(v_1\otimes v_2)$
is dominant.
Then for any $P\in U'{}(\Gb)$,
we have
$$P(v_1\otimes v_2)=v_1\otimes Pv_2+\cdots.$$
Hence the weights of
$\Usb v_2$ is contained in $\wt(v_2)+Q_-$.
Since $\Usb v_2=\Us v_2$ by Prop \ref{prop:bs},
$v_2$ is an extremal vector
with a dominant weight.
Similary so is $v_1$.
%
\QED

\newsection{Conjecture}\label{Conjecture}

We denote $\bigcup _{n>0}\Bbb C((q^{1/n}))$ by $\bk$
and $\bigcup _{n>0}\Bbb C[[q^{1/n}]]$  by $\bA$.
Hence $\bk$ is an algebraically closed field and
$\bA$ is a local ring.
For $a,\,b\in \bk^\times=\bk\setminus\{0\}$, we write $a\leq b$ if $a/b\in\bA$.

For $i\in I_0$, let $u_{i}$ denote the dominant extremal vector 
of $V(\varpi_{i})$.

\newcounter{new}
\setcounter{new}{\value{conj}}
\setcounter{conj}{\value{old}}
\addtocounter{conj}{-1}
\begin{conj}
Let $i_1,\ldots,i_l$ be elements of $I_0$ and
$a_{1},\ldots,a_{l}$ non-zero elements of $\bk$.
\hb
(1)\quad If $a_{1} \leq \cdots \leq a_{l}$, then 
$V(\varpi_{i_{1}})_{a_{1}}\otimes\cdots\otimes 
V(\varpi_{i_{l}})_{a_{l}}$ 
is generated by $u_{i_{1}}\otimes\cdots\otimes u_{i_{l}}$ 
as a $\Us_\bk$-module.
\hb
(2)\quad If $a_{1} \geq \cdots \geq a_{l}$,
then any non-zero $U'_q({\g g})_\bk$-submodule of 
$V(\varpi_{i_{1}})_{a_{1}}\otimes\cdots\otimes V(\varpi_{i_{l}})_{a_{l}}$ 
contains $u_{i_{1}}\otimes\cdots\otimes u_{i_{l}}$.
\end{conj}
\setcounter{conj}{\value{new}}

Note that (1) and (2) are dual statements and therefore they are equivalent.
One can compare (1) to the case of Verma modules and
(2) to the case of the dual of Verma modules.

Let us discuss several consequences of this conjecture.

For $i$,$j\in I_0$, there is an intertwiner
\eq
\Rn_{ij}(x,y):V(\varpi_i)_x\otimes V(\varpi_j)_y\to
V(\varpi_j)_y\otimes V(\varpi_i)_x
\endeq
We normalize this such that $R$ sends
$u_i\otimes u_j$ to $u_j\otimes u_i$.
Then we regard it as a rational function in $(x,y)$.
Since it is homogeneous,
its pole locus has the form $y/x=\on{constant}$.
We call it the {\em normalized $R$-matrix.}
By Corollary \ref{cor:irr} such an $\Rn_{ij}(x,y)$
is unique,

\begin{cor}\label{cor:pole}
If $a_{1} \leq a_{2}$,
the normalized $R$-matrix
$\Rn_{i,j}(x,y)$ does not have a pole
at $(x,y)=(a_1,a_2)$.
\end{cor}

\begin{pf}
Suppose that $\Rn_{i,j}(x,y)$ has a pole
at $(x,y)=(a_1,a_2)$.
Let $R'$ be the non-zero $\Us$-linear map
$V(\varpi_{i})_{a_{1}}\otimes V(\varpi_{j})_{a_{2}}
\to V(\varpi_{j})_{a_{2}}\otimes V(\varpi_{i})_{a_{1}}$
obtained after cancelling the poles of $\Rn_{i,j}(x,y)$. 
Then $R'(u_{i}\otimes u_{j}) = 0$,
and hence $\Im(R')$ does not have weight $\varpi_i+\varpi_j$.
On the other hand, Conjecture \ref{conj:main} (2)
implies that $\Im(R')$ contains $u_{j}\otimes u_{i}$, 
which is a contradiction. 
Hence $\Rn_{ij}(x,y)$ has no pole at $(a_{1},a_{2})$.
\end{pf}

\begin{cor}\label{cor:gen}
Let $K$ be a field extension of $k$, and
$i_1,\ldots,i_l\in I_0$,
$a_{1},\ldots,a_{l}\in K^\times=K\setminus\{0\}$.
\hb
$(1)$\quad Assume that $\Rn_{i_\nu,i_\mu}(x,y)$ does not have a pole
at $(x,y)=(a_\nu,a_\mu)$ for $1\le \nu<\mu\le l$. Then
$V(\varpi_{i_{1}})_{a_{1}}\otimes\cdots\otimes 
V(\varpi_{i_{l}})_{a_{l}}$ 
is generated by 
$u_{i_{1}}\otimes\cdots\otimes u_{i_{l}}$ 
as a $\Us_K$-module.
\hb
$(2)$\quad Assume that $\Rn_{i_\nu,i_\mu}(x,y)$ does not have a pole
at $(x,y)=(a_\nu,a_\mu)$ for $1\le \mu<\nu\le l$. Then
any non-zero $U'_q({\g g})_K$-submodule of 
$V(\varpi_{i_{1}})_{a_{1}}\otimes\cdots\otimes V(\varpi_{i_{l}})_{a_{l}}$ 
contains $u_{i_{1}}\otimes\cdots\otimes u_{i_{l}}$.
\end{cor}

\begin{pf}
We may assume that $K$ is generated by 
$a_{1},\ldots,a_{l}$ over $k$.
Since $\bk$ is an algebraically closed field with infinite transcendental
dimension over $k$,
there exists an embedding $K\hookrightarrow\bk$.
Hence we may assume $K=\bk$.

Since the proof of (2) is similar, we shall only prove
(1).
We prove (1) by induction on
the number of pairs $(\nu ,\mu)$
with $\nu < \mu$ and $a_{\nu}\not\le a_{\mu}$, which we denote by $n$.
If $n=0$, the assertion follows immediately from
Conjecture \ref{conj:main}.
If $n>0$, take
$\nu$ such that $a_{\nu}\not\le a_{\nu+1}$.
Hence $a_{\nu+1}\le a_{\nu}$.
Then Corollary \ref{cor:pole} implies that
$\Rn_{i_{\nu+1},i_{\nu}}(x,y)$ does not have a pole
at $(x,y)=(a_{\nu+1},a_\nu)$.
Since $\Rn_{i_{\nu},i_{\nu+1}}(x,y)$ does not have a pole
at $(x,y)=(a_{\nu},a_{\nu+1})$ by the assumption,
$$\Rn_{i_{\nu},i_{\nu+1}}(a_{\nu},a_{\nu+1}):
V(\varpi_{i_{\nu}})_{a_{\nu}}\otimes V(\varpi_{i_{\nu+1}})_{a_{\nu+1}}
\to
V(\varpi_{i_{\nu+1}})_{a_{\nu+1}}\otimes V(\varpi_{i_{\nu}})_{a_{\nu}}
$$
and
$$\Rn_{i_{\nu+1},i_{\nu}}(a_{\nu+1},a_{\nu}):
V(\varpi_{i_{\nu+1}})_{a_{\nu+1}}\otimes V(\varpi_{i_{\nu}})_{a_{\nu}}
\to
V(\varpi_{i_{\nu}})_{a_{\nu}}\otimes V(\varpi_{i_{\nu+1}})_{a_{\nu+1}}
$$
are inverse of each other.
Hence we can reduce the
original case to the case where $\nu$ and $\nu+1$ are exchanged,
in which $n$ is smaller than the original one by $1$.
Hence the induction proceeds.
\end{pf}

Assume the condition $(1)$ in the preceding Corollary \ref{cor:gen}.
Let
$R$ be the intertwiner
\eqn
&&R:V(\varpi_{i_{1}})_{a_{1}}\otimes\cdots\otimes V(\varpi_{i_{l}})_{a_{l}}
\to 
V(\varpi_{i_{l}})_{a_{l}}\otimes\cdots\otimes V(\varpi_{i_{1}})_{a_{1}}
\endeqn
sending
$u_{i_{1}}\otimes\cdots\otimes u_{i_{l}}$ to 
$u_{i_{l}}\otimes\cdots\otimes u_{i_{1}}$,
obtained as the product of
$\Rn_{i_\nu,i_\mu}(a_\nu,a_\mu)$ with $1\le \nu<\mu\le l$. 

\begin{cor}\label{cor:ir}
Under the condition $(1)$ in Corollary \ref{cor:gen},
$\Im(R)$ is irreducible.
\end{cor}
Note that the condition (1) is satisfied if $K=\bk$ and
$a_{1} \leq \cdots \leq a_{l}$, and hence we can apply the corollary.
\begin{pf}
By Corollary \ref{cor:gen} (1),
$\Im(R)$ is generated by the dominant extremal vector
$u_{i_{l}}\otimes\cdots\otimes u_{i_{1}}$.
Since any submodule of $\Im(R)$ contains
the same vector by Corollary \ref{cor:gen} (2),
$\Im(R)$ is irreducible.
\end{pf}

In fact, $\Im(R)$ is absolutely irreducible.
Let us recall that,
for a (not necessarily algebraically closed) field $K$ containing $k$,
a $\Us_K$-module $M$
finite-dimensional over $K$ is called absolutely irreducible
if the following equivalent conditions are satisfied.

\noindent
(1)\quad For some algebraically closed field $K'$
containing $K$, $K'\otimes_KM$ is an irreducible $\Us_{K'}$-module.
\smallskip
\hb
(2)\quad For any algebraically closed field $K'$
containing $K$, $K'\otimes_KM$ is an irreducible $\Us_{K'}$-module.
\smallskip
\hb
(3)\quad $M$ is irreducible and $\End_{\Us_K}(M)\cong K$.
\medskip
We denote by $\m$
the maximal ideal $\bigcup_{n>0}q^{1/n}\C[[q^{1/n}]]$
of $\bA$.

\Cor\label{Rn:pole}
For $a\in \bk^\times$,
$y/x=a$ is a pole of $\Rn_{ij}(x,y)$ if and only if
$a\in \m$ and $V(\varpi_i)\otimes V(\varpi_j)_a$ is reducible.
\encor
\proof
By Corollary \ref{cor:pole},
if $y/x=a$ is a pole of $\Rn_{ij}(x,y)$,
then $a\in\m$. By a similar argument to Corollary \ref{cor:pole},
the irreducibility of $V(\varpi_i)\otimes V(\varpi_j)_a$
implies that $y/x=a$ is not a pole of $\Rn_{ij}(x,y)$.
Now assume that $y/x=a\in\m$ is not a pole of $\Rn_{ij}(x,y)$.
Since $\Rn_{ji}(a,1)$ is well defined, $\Rn_{ij}(1,a)$ is invertible.
Hence $V(\varpi_i)\otimes V(\varpi_j)_a$ is irreducible
by Corollary \ref{cor:ir}.
\QED

\Cor\label{gen:irr}
Let $K$ be an algebraically closed field
containing $k$.
If $M$ and $M'$ are irreducible finite-dimensional integrable
$\Us_K$-modules,
then
$M\otimes M'_z$ is an irreducible $\Us_K$-module 
except finitely many $z\in K$.
\encor
\proof
Let $M$ (resp. $M'$) be the irreducible subquotient
of $V(\varpi_{i_1})_{a_1}\otimes\cdots\otimes V(\varpi_{i_m})_{a_m}$
(resp. $V(\varpi_{i'_1})_{a'_1}\otimes\cdots\otimes
V(\varpi_{i'_{m'}})_{a_{m'}}$)
such that
$\Rn_{i_\nu,i_\mu}(x,y)$
(resp. $\Rn_{i'_\nu,i'_\mu}(x,y)$)
does not have a pole
at $(x,y)=(a_\nu,a_\mu)$
(resp. $(x,y)=(a'_\nu,a'_\mu)$)
for $1\le\nu<\mu\le m$ (resp. $1\le\nu<\mu\le m'$).
Then Corollary \ref{cor:ir} implies that
$M$ is isomorphic to the image of the $R$-matrix
$$R:V\longrightarrow W,$$
where $V=V(\varpi_{i_1})_{a_1}\otimes\cdots\otimes V(\varpi_{i_m})_{a_m}$
and $W=V(\varpi_{i_m})_{a_m}\otimes\cdots\otimes V(\varpi_{i_1})_{a_1}$.
Similarly $M'$ is isomorphic to the image of
$$R': V'\longrightarrow W'\,,$$
where $V'=V(\varpi_{i'_1})_{a'_1}\otimes\cdots\otimes V(\varpi_{i'_{m'}})_{a'_{m'}}$
and
$W'=V(\varpi_{i'_{m'}})_{a'_{m'}}\otimes\cdots\otimes V(\varpi_{i'_1})_{a'_1}$.
If $z$ is generic , $\Rn_{i_\nu,i'_{\nu'}}(x,y)$
does not have a pole at $(x,y)=(a_\nu,za'_{\nu'})$
and $\Rn_{i'_{\nu'},i_\nu}(x,y)$
does not have a pole at $(x,y)=(za'_{\nu'},a_\nu)$.
Hence the $R$-matrix
$W\otimes W'_z\to W'_z\otimes W$ is an isomorphism.
Hence the image of the composition
$$
\begin{CD}
V\otimes V'_z@>{R\otimes R'_z}>>
W\otimes W'_z@>{\sim}>>W'_z\otimes W
\end{CD}
$$
is isomorphic to $M\otimes M'_z$ and
it is irreducible by 
Corollary \ref{cor:ir}.
\QED

Hence the intertwiner
$M\otimes M'_z\to M'_z\otimes M$ is unique up to constant.

We give a conjecture on the poles of the $R$-matrices.

\begin{conj}\label{conj:pole}
For $i,j\in I_0$, the pole of
the normalized $R$-matrix $\Rn_{ij}(x,y)$ has the form
$y/x=\pm q^n$ for $n\in\gamma^{-1}Z\ $ with $0<n\le (\delta,\rho)$
except $D^{(3)}_4$
(where $\gamma$ is defined in (\ref{def:gamma})).
In the $D^{(3)}_4$ case the third root of unity appears
in the coefficients.
\end{conj}

As seen in the appendix, this is true for $A^{(1)}_n$ and 
$C^{(1)}_n$.

We can also ask if the following statements are true.

\addtocounter{equation}{1}
\noindent
(\theequation)\quad 
$\Rn_{ij}(x,y)$ has only a simple pole.

\smallskip
\addtocounter{equation}{1}
\noindent
(\theequation)\quad If $(x,y)=(a,b)$ is a pole of $\Rn_{ij}(x,y)$,
then the kernel of 
$\Rn_{ji}(b,a):V(\varpi_j)_b\otimes V(\varpi_i)_a\to
V(\varpi_i)_a\otimes V(\varpi_j)_b$
is irreducible.

\newsection{Reduction of the conjecture}%
\label{Aux.Conj}

In this section we shall prove that 
Conjecture \ref{conj:main} follows from Conjecture \ref{conj2} below. 
Let $\m=\bigcup_{n>0}q^{1/n}\C[[q^{1/n}]]$ be the maximal ideal of $\bA$.

\begin{conj}\label{conj2}
For every $i \in I_{0}$, 
there exist $N \in {\Bbb N}, \ b_{1}, 
\cdots,b_{N},\ c_{1},\ldots,c_{N} 
\in \m\backslash \{ 0 \},\ s_{1},\ldots,s_{N},\ t_{1},\ldots,t_{N}
\in I_{0}$,
an irreducible finite-dimensional
$U'_q({\gtg})_\bk$-module $W_{\mu}$ and a $\Us_\bk$-linear map
${\varphi}_{{\mu}} : 
V(\varpi_{i})\otimes V(\varpi_{s_{\mu}})_{b_{\mu}} 
\longrightarrow 
V(\varpi _{t_{\mu}})_{c_{\mu}} \otimes W_{\mu}$ 
for $\mu$ with $1 \leq \mu \leq N$, 
satisfying the following conditions.
Define $F_{0} =\bigoplus_{\xi\not=-\varpi_{i^*}} V(\varpi _{i})_\xi$
(recall that
$-\varpi_{i^*}$ is the lowest weight vector of $V(\varpi _{i})$)
and 
$F_{{\mu}} = 
\{ v \in F_{{\mu}-1}|
{\varphi}_{{\mu}}(v \otimes u_{s_{\mu}}) = 0 \}$ 
for $0<\mu \leq N$.
\hb
(1)\quad $F_{{N}} = \bk u_{i}$.\hb
(2)\quad ${\varphi}_{{\mu}}(F_{{\mu}-1} \otimes u_{s_{\mu}}) 
\subset V(\varpi_{t_{\mu}})_{c_{\mu}} \otimes w _{\mu}$.\hb
(3)\quad $V(\varpi _{s_{\mu}})_{b_{\mu}}$ 
is not isomorphic to $V(\varpi _{t_{\mu}})_{c_{\mu}}$.\hb
(4)\quad $V(\varpi _{s_{\mu}})_{b_{\mu}}$ is not a component of $W_{\mu}$
(see Definition \ref{component}).\hb
Here $u_{s_{\mu}}$ and $w_{\mu}$ are dominant extremal vectors of 
$V(\varpi _{s_{\mu}})_{b_{\mu}}$ and $W_{\mu}$, respectively.
\end{conj}

\smallskip
Let us show that Conjecture \ref{conj2}
implies Conjecture \ref{conj:main} (2). 

%

For $a_{1},\ldots ,a_{p}\in\bk^{*}$,
let $P(a_{1},\ldots,a_{p})$ denote the following statement.
\smallskip
\begin{list}{$P(a_{1},\ldots,a_{p})$:}
\item
For indeterminates $x_{1},\ldots,x_{l}$,
any dominant extremal vector of the
$\Us_{\bk(x_{1},\ldots,x_{l})}$ -module
$V(\varpi _{j_{1}})_{x_{1}}\otimes\cdots\otimes V(\varpi _{j_{l}})_{x_{l}}
\otimes V(\varpi _{i_{1}})_{a_{1}}\otimes\cdots\otimes 
V(\varpi _{i_{p}})_{a_{p}}$ 
is a constant multiple of 
$u_{j_{1}}\otimes\cdots\otimes u_{j_{l}}\otimes 
u_{i_{1}}\otimes\cdots\otimes u_{i_{p}}$.
\end{list}
\smallskip
\noindent
Assuming Conjecture \ref{conj2},
we shall prove the following lemma.
\begin{lemma}\label{lemma:6}
If $a_{1},\ldots ,a_{p}\in \bk^{*}$ 
satisfy $a_{1}\geq\cdots\geq a_{p}$, 
then $P(a_{1},\ldots,a_{p})$ holds. 
\end{lemma}
Since any non-zero 
finite-dimensional module contains a dominant extremal vector, 
this lemma implies Conjecture \ref{conj:main} (2).
We shall prove this lemma by induction on $p$.
First assume $p\geq 1$.
Then $P(a_{1},\ldots,a_{p-1})$ holds
by the hypothesis of induction.
Set $K=\bk(x_{1},\ldots,x_{l})$.
Let $x$ be another indeterminate.
By the existence of $R$-matrix,
$V(\varpi_{j_{1}})_{x_{1}}\otimes\cdots\otimes
V(\varpi_{j_{l}})_{x_{l}}\otimes V(\varpi_{i_{1}})_{a_{1}}
\otimes\cdots\otimes V(\varpi_{i_{p-1}})_{a_{p-1}}
\otimes V(\varpi_{i_{p}})_{x}$ 
is isomorphic to
$V(\varpi_{j_{1}})_{x_{1}}\otimes\cdots\otimes
V(\varpi_{j_{l}})_{x_{l}}\otimes V(\varpi_{i_{p}})_{x}
\otimes V(\varpi_{i_{1}})_{a_{1}}
\otimes\cdots\otimes V(\varpi_{i_{p-1}})_{a_{p-1}}$.
Hence $P(a_{1},\ldots,a_{p-1})$ implies that
a dominant extremal vector of
the $\Us_{K(x)}$-module
$V(\varpi_{j_{1}})_{x_{1}}\otimes\cdots\otimes 
V(\varpi_{j_{l}})_{x_{l}}\otimes V(\varpi_{i_{1}})_{a_{1}}
\otimes\cdots\otimes V(\varpi_{i_{p-1}})_{a_{p-1}}
\otimes V(\varpi_{i_{p}})_{x}$ 
is a constant multiple of 
$u_{j_{1}}\otimes\cdots\otimes u_{j_{l}}\otimes u_{i_{1}}
\otimes\cdots\otimes u_{i_{p}}$. 
Then it follows that a dominant extremal vector of
$\Us_K$-module
$V(\varpi_{j_{1}})_{x_{1}}
\otimes\cdots\otimes V(\varpi_{j_{l}})_{x_{l}}
\otimes V(\varpi_{i_{1}})_{a_{1}}
\otimes\cdots\otimes 
V(\varpi_{i_{p-1}})_{a_{p-1}}
\otimes V(\varpi_{i_{p}})_{z}$ 
is a constant multiple of 
$u_{j_{1}}\otimes\cdots\otimes u_{j_{l}}
\otimes u_{i_{1}}\otimes\cdots\otimes u_{i_{p}}$ 
except for finitely many $z\in\bk$.
This means that $P(a_{1},\ldots,a_{p-1},z)$ 
holds except finitely many $z\in\bk$.
Arguing by induction on the order of the zero of $a_p$,
we may assume from the beginning
\eq\label{hyp:ind}
\mbox{$P(a_{1},\ldots,a_{p-1},z)$ holds
for any $z\in\m a_p\setminus\{0\}$.}&&
\endeq
Let $v$ be a dominant extremal vector 
of $\Us_K$-module
$V(\varpi_{j_{1}})_{x_{1}}\otimes\cdots\otimes 
V(\varpi_{j_{l}})_{x_{l}}\otimes V(\varpi_{i_{1}})_{a_{1}}
\otimes\cdots\otimes V(\varpi_{i_{p}})_{a_{p}}$.
We shall prove that $v$ is a constant multiple of 
$u_{j_{1}}\otimes\cdots\otimes u_{j_{l}}\otimes u_{i_{1}}
\otimes\cdots\otimes u_{i_{p}}$. 

We have $\varphi_0:V(\varpi_{i_p})\otimes V(\varpi_{{i_p}^*})_y\to \bk$
with $y=(-1)^{(\delta,\rho^\vee)}q^{(\delta,\rho)}$
by (\ref{duality}).
Set $V' = V({\varpi}_{j_{1}})_{x_{1}}
\otimes\cdots\otimes V({\varpi}_{j_{l}})_{x_{l}}
\otimes V({\varpi}_{i_{1}})_{a_{1}}
\otimes\cdots\otimes V({\varpi}_{i_{p-1}})_{a_{p-1}}$. 
Then we have a morphism
$$\id_{V'}\otimes(\varphi_0)_{a_p}:
V'\otimes V(\varpi_{i_p})_{a_p}\otimes V(\varpi_{{i_p}^*})_{a_py}
\to V'.$$
\begin{lemma}
We have $(\id_{V'} \otimes (\varphi _0)_{a_p})
(v \otimes u_{{i_p}^*}) = 0$.
\end{lemma}
\begin{pf}
Assume that
$w=(\id_{V'} \otimes (\varphi _0)_{a_p})
(v \otimes u_{{i_p}^*})\not= 0$.
Then $w$ is a dominant extremal vector of $V'$.
Hence $w$ is equal to
$u_{j_{1}}\otimes\cdots\otimes u_{j_{l}}\otimes u_{i_{1}}
\otimes\cdots\otimes u_{i_{p-1}}$
up to a constant multiple
by $P(a_1,\ldots,a_{p-1})$.
Therefore Theorem \ref{drin} implies that
$V(\varpi_{{i_p}^*})_{a_py}$ is isomorphic to
one of
$V({\varpi}_{j_{1}})_{x_{1}},\ldots,
V({\varpi}_{j_{l}})_{x_{l}},\ V({\varpi}_{i_{1}})_{a_{1}},
\cdots,V({\varpi}_{i_{p-1}})_{a_{p-1}}$.
This is a contradiction since $y\in qA$.
\end{pf}

Since $F_0=\{w\in V(\varpi_{i_p});\varphi_0(w\otimes u_{{i_p}^*})=0\}$,
we have $v\in V'\otimes (F_0)_{a_p}$.
Now we shall show $v\in  V'\otimes (F_\mu)_{a_p}$
by induction on $\mu$.
Applying Conjecture \ref{conj2} with $i=i_{p}$, 
we have $U'_q({\gtg})$-linear maps 
${\varphi}_{{\mu}} : 
V(\varpi_{i_{p}})\otimes 
V(\varpi_{s_{\mu}})_{b_{\mu}} 
\longrightarrow V(\varpi_{t_{\mu}})_{c_{\mu}} 
\otimes W_{\mu}$ for $1 \leq \mu \leq N$ 
satisfying the conditions $(1)-(4)$ in Conjecture \ref{conj2}.
Then this induces a homomorphism
$$\id_{V'} \otimes(\varphi _{\mu})_{a_p}
:V'\otimes V(\varpi_{i_{p}})_{a_p}\otimes
V(\varpi_{s_{\mu}})_{a_pb_{\mu}}
\longrightarrow V'\otimes V(\varpi_{t_{\mu}})_{a_pc_{\mu}}
\otimes (W_{\mu})_{a_p}\,.$$
Suppose that $v \in V' \otimes (F_{{\mu}-1})_{a_p}$, 
which is the case when $\mu = 1$.

\begin{lemma}
We have $(\id_{V'} \otimes 
(\varphi _{\mu})_{a_p})(v \otimes u_{s_{\mu}}) = 0$.
\end{lemma}

\begin{pf}
The proof is similar to the one of the preceding lemma.
Suppose that $w=(\id_{V'} \otimes
(\varphi _{\mu})_{a_p})(v \otimes u_{s_{\mu}})$ is not zero.
Write $w$ as $v'' \otimes w_{\mu}$ in virtue of the condition (2)
in Conjecture \ref{conj2}, 
where $v''$ is a non-zero vector of $V' \otimes 
V({\varpi}_{t_{\mu}})_{a_pc_{\mu}}$.
Since $v \otimes u_{s_{\mu}}$ is extremal, 
so is $v'' \otimes w_{\mu}$. 
Hence $v''$ is a dominant extremal vector by Lemma \ref{ext:mod}. 
Since $a_pc_{\mu}\in\m a_p$,
the property
$P(a_1,\ldots,a_{p-1},a_pc_\mu)$ holds by (\ref{hyp:ind}),
and hence $v''$ is a nonzero scalar multiple of 
$u_{j_{1}}\otimes\cdots\otimes u_{j_{l}}\otimes u_{i_{1}}
\otimes\cdots\otimes u_{i_{p}}\otimes u_{t_{\mu}}$.
Then Theorem \ref{drin} implies that $V({\varpi}_{s_{\mu}})_{a_pb_{\mu}}$ 
is isomorphic to one of
$V({\varpi}_{j_{1}})_{x_{1}},\ldots,$$
V({\varpi}_{j_{l}})_{x_{l}}$, 
$V({\varpi}_{i_{1}})_{a_{1}},\cdots,$
$V({\varpi}_{i_{p-1}})_{a_{p-1}}$, 
$V({\varpi}_{t_{\mu}})_{a_pc_{\mu}}$
or to a component of $ (W_{\mu})_{a_p}$.
It, however, is not the case because of the conditions
(3), (4) in Conjecture \ref{conj2} and $a_pb_{\mu}\in \m a_{p}$.
\end{pf}

By this we have $v \in V' \otimes (F_{\mu})_{a_p}$. 
Applying this process successively, 
we obtain $v \in V' \otimes (F_{N})_{a_p}$. 
Hence we have $v \in V' \otimes u_{i_{p}}$ 
by the condition (1) in Conjecture \ref{conj2}. 
Write $v$ as $v' \otimes u_{i_{p}}$, 
where $v'$ is a nonzero vector of $V'$. 
Lemma \ref{ext:mod} implies that $v'$ is dominant and extremal.
Therefore $v'$ is a nonzero scalar multiple of 
$u_{j_{1}}\otimes\cdots\otimes u_{j_{l}}\otimes 
u_{i_{1}}\otimes\cdots\otimes u_{i_{p-1}}$ 
by the induction hypothesis on $p$. 
We have deduced the $p$ case from the $p-1$ case.

It remains to prove $p=0$ case, which follows from the following lemma.

\begin{lemma}
Any dominant extremal vector of the $\Us_{\bk(x_{1},\ldots,x_{l})}$-module
$V(\varpi_{i_{1}})_{x_{1}}\otimes\cdots\otimes 
V(\varpi_{i_{l}})_{x_{l}}$ is a constant multiple of 
$u_{i_{1}}\otimes\cdots\otimes u_{i_{l}}$.
Here $x_{1},\ldots,x_{l}$ are indeterminates.
\end{lemma}

\begin{pf}
It is enough to prove the assertion with $x_{1}=\cdots =x_{l}=1$.
Let $V$ denote $V(\varpi_{i_{1}})\otimes\cdots\otimes V(\varpi_{i_{l}})$. 
By Corollary \ref{cor:irr}, 
$V$ is irreducible.
Suppose now that $V$ has a dominant extremal vector $v$
that is not a constant multiple of
$u_{i_{1}}\otimes\cdots\otimes u_{i_{l}}$. 
Then $U'_q({\gtg})v$ does not contain 
$u_{i_{1}}\otimes\cdots\otimes u_{i_{l}}$ 
since $\wt(U'_q({\gtg})v)\subset 
\wt(v) + \sum _{i \in I_{0}}\Z_{\leq 0}\cl(\alpha _{i})$, 
which is a contradiction. 
\end{pf}

Thus we have proved 
\Prop Conjecture \ref{conj2} implies 
Conjecture \ref{conj:main}.
\enprop

\newsection{Proof of Conjecture \ref{conj2} for $A^{(1)}_n$ and $C^{(1)}_n$}%
\label{AC}
In this section, we shall prove the following theorem.
\begin{thm}
Conjecture \ref{conj2} holds if $\gtg$ is $A^{(1)}_{n}$ or $C^{(1)}_{n}$.
\end{thm}

\subsection{$A^{(1)}_{n-1}$ Case}
For the fundamental representations of $\uasl$,
see Appendix \ref{A^1_n}.
We identify crystal bases of the fundamental representations 
with the corresponding global bases.


Let us prove Conjecture \ref{conj2}.

Since the $i=n-1$ case can be 
reduced to the case $i=1$ by the Dynkin diagram automorphism,
we assume $1 \leq i<n-1$.
Set $N = i$. 
For $1 \leq\mu\leq N=i$, take 
$s_{\mu}= \mu,\ t_{\mu}= i + 1,\ b_{\mu}=(-q)^{i- \mu +2},
\ c_{\mu}= -q,\ W_{\mu} = V({\varpi}_{\mu -1})_{(-q)^{i- \mu + 1}}$ 
and define 
${\varphi}_{i,{\mu}} : 
V(\varpi_{i})\otimes V(\varpi_{\mu})_{(-q)^{i- \mu + 2}} 
\longrightarrow 
V(\varpi_{i+1})_{-q}\otimes W_{\mu}$ 
as the composition (see Lemma \ref{A1}):
\begin{equation}
  \begin{CD}
      V(\varpi_{i})\otimes V(\varpi_{\mu})_{(-q)^{i- \mu + 2}}  
@>{V(\varpi _{i})}\otimes (i_{1,\mu-1})_{(-q)^{i-\mu +2}}>>
     V(\varpi_{i})\otimes V(\varpi_{1})_{(-q)^{i+1}} \otimes W_{\mu}   \\
     @.       @VV{{(p_{i,1})_{-q}}\otimes {W_{\mu}}}V \\
       @.     V(\varpi_{i+1})_{-q}\otimes W_{\mu}.
 \end{CD}
\end{equation}
%
%
Then it is easy to check that Conjecture \ref{conj2} holds with 
\beqn
F_{\mu}&=& \bigoplus_{\mu<a_{\mu+1}<\cdots a_i\le n}
 k(1,\ldots ,\mu ,\ a_{\mu +1},\ldots ,a_{i})\,.\\
\endeqn

\subsection{$C^{(1)}_{n}$ Case}
For the fundamental representations of
$U_q'(C^{(1)}_{n})$, see Appendix \ref{C_n^{(1)}}.

%
For $1\le i<n$, let
$p_i:V(\varpi_i)\otimes V(\varpi_1)_{(-q_s)^{i+1}}\to V(\varpi_{i+1})_{-q_s}$
be $(p_{i,1})_{-q_s}$.
Let
$p_n:V(\varpi_n)\otimes V(\varpi_1)_{(-q_s)^{n+3}}\to V(\varpi_{n-1})_{-q_s}$
be the composition
\eqn
&&\begin{CD}
V(\varpi_n)\otimes V(\varpi_1)_{(-q_s)^{n+3}}@.{}\\
@V{i_{n-1,1}\otimes{V(\varpi_1)_{(-q_s)^{n+3}}}}VV\\
V(\varpi_{n-1})_{-q_s}\otimes V(\varpi_1)_{(-q_s)^{1-n}}\otimes 
V(\varpi_1)_{(-q_s)^{n+3}}@>{{V(\varpi_{n-1})_{-q_s}}\otimes \tr}
>>V(\varpi_{n-1})_{-q_s}.
\end{CD}
\endeqn
Here $\tr$ is given in (\ref{C:dual}).

For $1\leq i\leq n-1$, set $N=i$. 
For $ 1\leq\mu\leq N=i$, we set 
$s_{\mu}=\mu$, $t_{\mu}= i + 1$, 
$b_{\mu}= (-q_s)^{i- \mu +2}$, $c_{\mu}= -q_s$
and $W_{\mu}=V({\varpi}_{\mu -1})_{(-q_s)^{i- \mu + 1}}$.
We define 
${\varphi}_{i,{\mu}}:V(\varpi _{i})\otimes 
V(\varpi _{\mu})_{(-q_s)^{i- \mu +2}} 
\longrightarrow V(\varpi _{i+1})_{-q_s} \otimes W_{\mu}$ 
as the composition:

\begin{equation}
\begin{CD}
V(\varpi _{i})\otimes V(\varpi _
{\mu})_{(-q_s)^{i- \mu +2}} @>{{V(\varpi_{i})}\otimes
(i_{1,\mu-1})_{(-q_s)^{i-\mu +2}}}>>   V(\varpi _{i})\otimes 
V(\varpi _{1})_{(-q_s)^{i+1}}\otimes W_{\mu} \\
  @.  @VV{p_{i}\otimes{W_{\mu}}}V  \\
@.   V(\varpi _{i+1})_{-q_s} \otimes W_{\mu}. 
\end{CD}
\end{equation}
Note that $b_{\mu},\ c_{\mu}\in q_sA$. \par 

For $i=n$, set $N=n$. For $1\leq\mu\leq n$, 
we set $s_{\mu}= \mu$, $t_{\mu}= n-1$,
$b_{\mu}=(-q_s)^{n- \mu +4}$, $c_{\mu}= -q_s$ and 
$W_{\mu}=V({\varpi}_{\mu -1})_{(-q_s)^{n- \mu + 3}}$.
We define 
${\varphi}_{n,{\mu}}:V(\varpi _{n})\otimes 
V(\varpi_{\mu})_{(-q_s)^{n-\mu +4}}$
$ \longrightarrow V(\varpi _{i+1})_{-q_s} \otimes W_{\mu}$ as the composition ;
 
\begin{equation}
\begin{CD}
V(\varpi _{n})\otimes V(\varpi_{\mu})_{(-q_s)^{n-\mu +4}} 
@>{{V(\varpi_{n})}\otimes (i_{1,\mu-1})_{(-q_s)^{n-\mu+4}}}>> 
  V(\varpi _{n})\otimes V(\varpi _{1})_{(-q_s)^{n+3}} 
\otimes W_{\mu}   \\
    @.   @VV{p_{n}\otimes{W_{\mu}}}V   \\
      @.     V(\varpi _{n-1})_{-q_s} \otimes W_{\mu}. 
\end{CD}
\end{equation}
Note that $b_{\mu},\ c_{\mu}\in q_sA$.

Then we have
$$F_\mu=\{v\in V(\varpi_i);
p_i(v\otimes G(j))=0\quad\mbox{for $1\le j\le \mu$}\}.
$$
Then Conjecture \ref{conj2} easily follows from
the following lemma.

\begin{lemma}\label{lemma:14}
Fix $1\leq i\leq n$. 
Then 
\eq\label{eq:C}
&&\{ v\in V({\varpi}_{i})\ | \ p_{i}(v \otimes G(j)) = 0
\quad\text{for all $1 \leq j \leq i$} \} = kG(1,\ldots ,i).
\endeq
\end{lemma}
\begin{pf}
Let $E$ be the left-hand-side of (\ref{eq:C}).
Then $E$ is invariant 
by $e_{k}$ for any $k\in I_0$. 
Let us prove $E_\lam=0$ by induction on 
the weight $\lam\not=\varpi_i$.
We can easily check the assertion when 
$\lam=\varpi_i-\alpha_i$,
since $V(\varpi_i)_\lam=kf_iu_i$.
If a weight $\lam$ of $V(\varpi_i)$ is not $\varpi_i-\alpha_i$, 
then $\lam+\alpha_k\not=\varpi_i$ for any $k\in I_0$.
Therefore any $v\in E_\lam$ satisfies
$e_{k}v=0$ for all $k\in I_0$ by the induction hypothesis. 
This implies $v=0$.
\end{pf}

\appendix
\newsection{Universal R-matrix}
In this appendix we shall calculate the normalized and 
universal $R$-matrices of 
$\Us$ for the fundamental representations following 
a variant of the recipe
of Frenkel-Reshetikhin \cite{FR}
in the $A^{(1)}_{n-1}$ and $C^{(1)}_{n}$ cases.

Let us choose the following universal $R$-matrix.
Let us take a base $P_\nu$ of $U_q^+(\gtg)$
and $Q_\nu$ of $U_q^-(\gtg)$
dual to each other with respect a suitable coupling between
$U_q^+(\gtg)$ and $U_q^-(\gtg)$.
Then for $\Us$-modules $M$ and $N$
define
\eq\label{def:univ}
R_{MN}^\univ(u\otimes v)=q^{(\wt(u),\wt(v))}
\sum_\nu P_\nu v\otimes Q_\nu u\,,
\endeq
so that
$R_{MN}^\univ$ gives a $\Us$-linear homomorphism
from $M\otimes N$ to $N\otimes M$
provided the infinite sum has a meaning.
If $M$ and $N$ are finite-dimensional integrable modules,
then $R_{M,N_z}^\univ$
converges in the $z$-adic topology.
%
%
%
%
%
%
The existence of the universal $R$-matrix for $(M,N)$ is proved by 
\cite{Drin2} (see also \cite{T}).
For a scalar $a$, the composition
$$(R_{M,N}^\univ)_a
:M_a\otimes N_a\cong (M\otimes N)_a\to
(N\otimes M)_a\cong M_a\otimes N_a
$$
is equal to
$R_{M_a,N_a}^\univ$, and we sometimes confuse them.

For irreducible $\Us$-modules $M$ and $N$, let us denote
by $\Rn_{MN}(z)$ the $R$ matrix
$M\otimes N_z\to N_z\otimes M$
normalized by $\Rn_{MN}(z)(u\otimes v)=v\otimes u$
for dominant extremal vectors $u$ (resp. $v$)
of $M$ (resp. $N$).
Let $d_{MN}(z)$ be
a denominator of $\Rn_{MN}(z)$.
Namely $c(z)\in k[z,z^{-1}]$
is divisible by $d_{MN}(z)$ if and only if
$c(z)\Rn_{MN}(z)$ has no poles.
Then $d_{MN}(z)$ is uniquely determined modulo
$k[z,z^{-1}]^\times$.
Here $k[z,z^{-1}]^\times$ is the set of invertible elements of
$k[z,z^{-1}]$.
Hence
\eq
&&k[z,z^{-1}]^\times=\{cz^n;n\in\Z,\ c\in k\setminus\{0\}\}.
\endeq 
Since the intertwiner from $M\otimes N_z$ to $N_z\otimes M$
is unique up to a constant multiple by Corollary \ref{gen:irr},
we can write
\eq
\Ru_{MN}(z)=a_{MN}(z)\Rn_{MN}(z).
\endeq
If $\lam$ and $\mu$ are the dominant extremal weight of $M$ and $N$ 
respectively,
we have
\eq\label{con:init}
a_{MN}(z)\in q^{(\lam,\mu)}(1+zk[[z]]).
\endeq

For $i,j\in I_0$, we denote 
$\Ru_{ij}(z)=\Ru_{V(\varpi_i)V(\varpi_j)}(z)$,
$\Rn_{ij}(z)=\Rn_{V(\varpi_i)V(\varpi_j)}(z)$,
$a_{ij}(z)=a_{V(\varpi_i)V(\varpi_j)}(z)$ and
$d_{ij}(z)=d_{V(\varpi_i)V(\varpi_j)}(z)$.

For a finite-dimensional $\Us$-module $M$,
let $M^*$ be the left dual of $M$ and
${}^*M$ the right dual of $M$.
Hence we have
\eqn
\begin{CD}
M^*\otimes M@>{\tr}>>k\\
M\otimes {}^*M@>{\tr}>>k
\end{CD}
&\qquad&
\begin{CD}
k@>{\iota}>>M\otimes M^*\\
k@>{\iota}>>{}^*M\otimes M.
\end{CD}
\endeqn
We have
\eq\label{bidual}
&&M^*{}^*\cong M_{q^{-2(\delta,\rho)}}\quad\text{and}\quad
{}^*{}^*M\cong M_{q^{2(\delta,\rho)}}.
\endeq
We have
\eqn
&&V(\varpi_i)^*\cong V(\varpi_{i^*})_{\p^{-1}}\quad\text{and}\quad
{}^*V(\varpi_i)\cong V(\varpi_{i^*})_{\p},
\endeqn
where $\p=(-1)^{(c,\rho^\vee)}q^{(\delta,\rho)}$.

Let $a\mapsto\overline{a}$ be the ring automorphism
of $\Us$ given by $\overline{q}=q^{-1}$,
$(e_i)^-=e_i$, $(f_i)^-=f_i$, $q(h)^-=q(-h)$.
For a $\Us$-module $M$, let $M^-$ be the $\Us$-module
whose underlying vector space
is $M$ with the new action $\Us\overset{-}{\To}\Us\To \End(M)$.
Then $(M\otimes N)^-\cong N^-\otimes M^-$ and
$V(\varpi_i)^-\cong V(\varpi_i)$.
Hence we have
\eq\label{reverse}
d_{ji}(z)\equiv d_{ij}(z^{-1})^-\quad\mod k[z,z^{-1}]^\times.
\endeq
The conjecture \ref{conj:pole} implies
\eq
d_{ji}(z)\equiv d_{ij}(z)\quad\mod k[z,z^{-1}]^\times.
\endeq

\Prop
For irreducible finite-dimensional integrable $\Us$-modules
$V$ and $W$, we have
\eq\label{eq0:univ}
&&a_{V,W}(z)a_{{}^*V,W}(z)
\equiv \dfrac{d_{VW}(z)}{d_{W,{}^*V}(z^{-1})}
\quad\mod k[z,z^{-1}]^\times.
\endeq
\enprop

\begin{pf}
For a $\Us$-linear homomorphism
$\phi:V\otimes W_z\to W_z\otimes V$,
we shall define $\Tr(\phi):W_z\otimes {}^*V\to {}^*V\otimes W_z$
as the composition
\eqn
\begin{CD}
W_z\otimes {}^*V@>{\iota\otimes W_z\otimes {}^*V}>>
{}^*V\otimes V\otimes W_z\otimes {}^*V
@>{{}^*V\otimes \phi\otimes {}^*V}>>
{}^*V\otimes W_z\otimes V\otimes {}^*V
\end{CD}&&\\
&&\hspace{-100pt}
\begin{CD}
@>{{}^*V\otimes W_z\otimes\tr}>>
{}^*V\otimes W_z.
\end{CD}
\endeqn
The correspondence
$\phi\mapsto \Tr(\phi)$
gives an isomorphism
\eq
\Hom(V\otimes W_z, W_z\otimes V)\overset{\sim}{\longrightarrow}
\Hom(W_z\otimes {}^*V,{}^*V\otimes W_z).
\endeq
If we consider them as modules over $k[z,z^{-1}]$, then
$\Hom(V\otimes W_z, W_z\otimes V)$ is generated by
$d_{VW}(z)\Rn_{VW}(z)$, and
$\Hom(W_z\otimes {}^*V,{}^*V\otimes W_z)$ is generated by
$d_{W,{}^*V}(z^{-1})\Rn_{W,{}^*V}(z^{-1})$.
Hence we have
\eqn
\Tr(d_{VW}(z)\Rn_{VW}(z))\equiv
d_{W,{}^*V}(z^{-1})\Rn_{W,{}^*V}(z^{-1})
\qquad\mod k[z,z^{-1}]^\times.
\endeqn
Then the result follows from
$\Rn_{W,{}^*V}(z^{-1})
=(\Rn_{{}^*V,W}(z))^{-1}$
and a well known result
$\Tr(\Ru_{VW}(z))
=(\Ru_{{}^*V,W}(z))^{-1}$
(see \cite{FR}).
\end{pf}
This proposition implies
\eq\label{eq:univ}
&&{a_{i,j}(z)}{a_{i^*,j}({\p}^{-1}z)}
\equiv \dfrac{d_{i,j}(z)}{d_{j,i^*}(\p z^{-1})}
\quad\mod k[z,z^{-1}]^\times.
\endeq
Applying (\ref{eq0:univ}) with ${}^*V$ instead of $V$, we have
\eq
&&a_{{}^*V,W}(z)a_{{}^{**}V,W}(z)
\equiv \dfrac{d_{{}^*V,W}(z)}{d_{W,{}^{**}V}(z^{-1})}
\quad\mod k[z,z^{-1}]^\times.
\endeq
Using (\ref{bidual}) we obtain the $q$-difference equation
\eq\label{eq:diff}
&&\dfrac{a_{VW}(z)}{a_{VW}(q^{-2(\delta,\rho)}z)}
\equiv
\dfrac{d_{VW}(z)d_{WV}(q^{2(\delta,\rho)}z^{-1})}
{d_{W,{}^*V}(z^{-1})d_{{}^*V,W}(z)}
\quad\mod k[z,z^{-1}]^\times.
\endeq
Write
\eqn
&&d_{ji}(z)=\prod_\nu(z-x_\nu)\quad\mbox{and}\quad
d_{j,i^*}(z)=\prod_\nu(z-y_\nu).
\endeqn
Then by (\ref{reverse}), we have
\eqn
&&d_{ij}(z)=\prod_\nu(z-\ol{x_\nu}^{-1})\quad\mbox{and}\quad
d_{i^*,j}(z)=\prod_\nu(z-\ol{y_\nu}^{-1}).
\endeqn
Then using (\ref{con:init}), 
we can solve the $q$-difference equation
(\ref{eq:diff}),
\eq\label{gen:a}
&&a_{ij}(z)=q^{(\varpi_i,\varpi_j)}
\dfrac{\prod_\nu(\p y_\nu z;\p^2)_\infty(\p\overline{y_\nu}z;\p^2)_\infty}
{\prod_\nu(x_\nu z;\p^2)_\infty(\p^2\overline{x_\nu}z;\p^2)_\infty}\ .
\endeq
Here $\p=(-1)^{(\delta,\rho^\vee)}q^{(\delta,\rho)}$
and $(z;q)_\infty=\prod_{n=0}^\infty(1-q^nz)$.

We are going to determine $d_{ij}(z)$ and $a_{ij}(z)$
in the $A^{(1)}_{n-1}$ and $C^{(1)}_{n}$ cases.

\begin{rem}
We can see easily
\eq
&&d_{V^*,W^*}(z)\equiv d_{{}^*V,{}^*W}(z)\equiv d_{V,W}(z).
\endeq
Hence
\eq
&&a_{V^*,W^*}(z)=a_{{}^*V,{}^*W}(z)=a_{V,W}(z),
\endeq
and
\eq
&&d_{i^*,j^*}(z)\equiv d_{i,j}(z).
\endeq
\end{rem}

\section{$A_{n-1}^{(1)}$ case}
 We shall review the fundamental representations
and $R$-matrices for $A_{n-1}^{(1)}$ .
\subsection{Fundamental representations}\label{A^1_n}
The root data of $\Gg=A_{n-1}^{(1)}$ are as follows.
\eqn
&&I=\{0,1,\ldots,n-1\}\\
&&(\alpha_i,\alpha_j)=
\begin{cases}
2&\mbox{if $i=j$}\\
-\delta(i\equiv j+1\, \mod\ n)-\delta(i\equiv j-1\, \mod\ n)
&\mbox{otherwise}
\end{cases}\\
&&\delta=\alpha_0+\cdots+\alpha_{n-1}\,,\\
&&c=h_0+\cdots+h_{n-1}\,,\\
&&(\delta,\rho)=(\delta,\rho^\vee)=n.
\endeqn
Here for the statement $P$, we define  $\delta(P)=1$ or $-1$
according that $P$ is true or false.

Hence by (\ref{duality}) the duality morphisms are given by
\eqn
k\overset{\iota}{\To}V(\varpi_{n-i})_{(-q)^{n}}\otimes V(\varpi_i)
\quad\text{and}\quad
V(\varpi_i)\otimes V(\varpi_{n-i})_{(-q)^{n}}\overset{\tr}{\To}k.
\endeqn

By \cite{(KMN)^2_2}, the vectors of the crystal base $B_{k}$
of the fundamental representation
$V(\varpi _{k})$ $(1\leq k\leq n-1)$ are labeled by
the subsets of $\Z/n\Z=\{1,\ldots,n\}$ with exactly $k$ elements.
For $0\le i\le n-1$ and $K\subset \Z/n\Z$, we have
\eqn
\te_i(K)&=&
\begin{cases}
(K\setminus\{i+1\})\cup\{i\}
&\mbox{if $i+1\in K$ and $i\not\in K$,}\\
0&\mbox{otherwise,}
\end{cases}\\
\tf_i(K)&=&
\begin{cases}
(K\setminus\{i\})\cup\{i+1\}
&\mbox{if $i\in K$ and $i+1\not\in K$,}\\
0&\mbox{otherwise.}
\end{cases}
\endeqn

In the case of the fundamental representations of $\uasl$,
all the weights are extremal.
Therefore we have $e_iG(b)=G(\te_ib)$ and 
$f_iG(b)=G(\tf_ib)$ for every $b$ in the crystal base.
Here $G(b)$ is the corresponding global base.
Hence we can and do
identify its crystal bases with the corresponding global bases.

We have
\eqn
&&t_iK=q^{\delta(i\in K)-\delta(i+1\in K)}K\,.
\endeqn

We present a lemma that is easily verified by calculation.
\begin{lemma}\label{A1}
For $j,k\ge 0$ such that $j+k\le n$,
there exist following non-zero $\uasl$-linear 
homomorphisms.
\begin{enumerate}
\item
$i_{j,k}:V(\varpi_{j+k})\longrightarrow
V(\varpi_j)_{(-q)^{k}}\otimes V(\varpi_k)_{(-q)^{-j}}$
\hb
given by
\eqn
&&i_{j,k}(M)=\sum_{\sharp J=j,\sharp K=k \atop M=J\cup K,\,J\cap K=\emptyset}
(-q)^{\psi(J,K)}J\otimes K\,.
\endeqn
Here $\psi(J,K)=\sharp\{(\nu,\mu)\in J\times K\,;\,\nu>\mu\}$.
\item
$p_{j,k}:V(\varpi_j)_{(-q)^{-k}}\otimes V(\varpi_k)_{(-q)^{j}}
\longrightarrow V(\varpi_{j+k})$
\hb
given by
\eqn
&&p_{j,k}(J\otimes K)=
\begin{cases}
(-q)^{\psi(J,K)}(J\cup K)&\mbox{if $J\cap K=\emptyset$}\\
0&\mbox{if $J\cap K\not=\emptyset$.}
\end{cases}
\endeqn
\end{enumerate}
Here $V(\varpi_0)$ and $V(\varpi_n)$ are 
understood to be the trivial representation.
\end{lemma}

\subsection{$R$-matrices}
We shall recall the result of Date-Okado \cite{DO}.
\begin{prop}[\cite{DO}]
For $k,l\in I_0$
\eq\label{A:d}
&d_{kl}(z)=\mathop{\displaystyle\prod}\limits_{\nu=1} ^{\min(k,l,n-k,n-l)}
(z-(-q)^{2\nu+|k-l|}).
\endeq
\enprop
%

The universal $R$-matrices can be easily obtained by
(\ref{gen:a}) and (\ref{A:d}).

\begin{prop}[\cite{DO}]For $k,l\in I_0=\{1,\ldots,n-1\}$, we have
\beqn
a_{kl}(z)=
q^{\min(k,l)-kl/n}\dfrac{((-q)^{|k-l|}z;q^{2n})_\infty
((-q)^{2n-|k-l|}z; q^{2n})_\infty}
{((-q)^{k+l}z;q^{2n})_\infty((-q)^{2n-k-l}z; q^{2n})_\infty}\,.
\endeqn
\end{prop}

\section{$C_{n}^{(1)}$ case}

\subsection{Fundamental representations}\label{C_n^{(1)}}
The Dynkin diagram of $C_{n}^{(1)}$ is
\newcommand{\toru}{\hspace{-10pt}}
$$
\begin{matrix}
{0}&&1&&2&&n-1&&n\\
{\bigcirc}&\toru\Longrightarrow\toru&\bigcirc&
\toru\rule[2.5pt]{20pt}{.5pt}\toru
&\bigcirc&\toru\rule[2.5pt]{20pt}{.5pt}\cdots\cdots\rule[2.5pt]{20pt}{.5pt}
\hspace{-15pt}&\bigcirc&
\hspace{-20pt}\Longleftarrow\toru&\bigcirc\\
-2\eps_1&&\eps_1-\eps_2&&\eps_2-\eps_3&&\eps_{n-1}-\eps_n&&2\eps_n
\end{matrix}
$$
Here $(\eps_i)_{i=1,\ldots,n}$ is an orthogonal basis
of $\clwt$ such that $(\eps_i,\eps_i)=1/2$.
We have
\eqn
q_i&=&
\begin{cases}
q&\mbox{if $i=0$ or $n$}\\
q^{1/2}&\mbox{if $1\le i<n$,}
\end{cases}\\
\delta&=&\alpha_0+2(\alpha_1+\cdots\alpha_{n-1})+\alpha_n\,,\\
c&=&h_0+h_1+\cdots+h_n\,,\\
(\delta,\rho)&=&n+1\,,\\
\lan\rho^\vee,\delta\ran&=&2n\,,\\
\varpi_i&=&\Lambda_i-\Lambda_0=\eps_1+\cdots+\eps_i\,.
\endeqn
We set $q_s=q^{1/2}$.
Hence by (\ref{duality}) the duality morphisms are given by
\eq\label{C:dual}
&&k\overset{\iota}{\To}V(\varpi_i)_{q_s^{2(n+1)}}\otimes V(\varpi_i)
\quad\text{and}\quad
V(\varpi_i)\otimes V(\varpi_i)_{q_s^{2(n+1)}}\overset{\tr}{\To}k.
\endeq
We review the crystal base $(L_{k},\ B_{k})$ 
of the fundamental representation 
$V(\varpi _{k})$ $(1\leq k\leq n)$ 
of $U'_q({C^{(1)}_{n}})$. 
Recall that  $V(\varpi _{k})$ is as 
a $U_q(C_{n})$-module isomorphic to 
the $k$-th fundamental representation of $U_q(C_{n})$. 
Hence by \cite{KN}, $B_{k}$ is labeled by 
$$
\aligned
\{ (m_{i})_{i=1}^{k}\ |\ m_{1} \prec \cdots \prec m_{k},
&\ m_{i}\in \{1,\ldots,n,\ \overset{-}{n},\ldots,\overset{-}{1}\}, \\ 
  &i+(k-j+1)\leq m_{i}\quad\text{if}\quad m_{i}=\overset{-}{m_{j}}\ (i<j)\},
\endaligned
$$
where the ordering on 
$\{1,\ldots,n,\ \ol{n},\ldots,\ol{1}\}$ 
is defined by 
\begin{equation}
1\prec\cdots \prec n \prec\ol{n} \prec\cdots \prec\ol{1}.
\end{equation}
On $B_{k}$ the actions of $\tf_{i}$ and 
$\te_{i}$ with $0 \leq i \leq n$ are defined as follows.
As for $i \neq 0$, write $i$ and $\ol{i+1}$ as $+$, $i+1$ 
and $\ol{i}$ as $-$, and others as $0$. 
Then first ignore $0$ and next ignore $+-$.
Then $\tf_{i}b$ is obtained 
by replacing the leftmost $+$ with $-$ and 
$\te_{i}b$ is obtained by replacing the rightmost $-$ 
with $+$.

\Lemma
If $b$ is of the form $(1,\ a_{1},\ldots ,a_{k-1})$, 
then $\te_{0}b=(a_{1},\ldots ,a_{k-1},\ \ol{1})$. 
Otherwise $\te_{0}b=0$.\hb
If $b$ is of the form $(a_{1},\ldots,a_{k-1},\ \ol{1})$, 
then $\tf_{0}b=(1,\ a_{1},\ldots ,a_{k-1})$.  Otherwise $\tf_{0}b=0$.
\enlemma
\proof
It is easy to check that
$B_k$ is a regular crystal
with this definition of 
$\te_{0}$ and $\tf_{0}$.
Set $J=\{1,2,\ldots,n-1\}\subset I$. Then $B_k$ decomposes,
as a crystal over $\Gg_J\simeq A_{n-1}$,
into irreducible components 
with multiplicity $1$.
Hence there is a unique way to draw $0$-arrows on 
the crystal $B_k$ over $C_n$.
\QED

The following proposition can be checked by a direct calculation.

\Prop
For $\mu,\nu$ with $\mu+\nu\le n$,
there exist following non-zero $U'_q({C^{(1)}_{n}})$-linear maps:
\eqn
&&i_{\mu\nu} : V({\varpi}_{\mu+\nu})
\longrightarrow V({\varpi}_{\mu})_{(-q_s)^{\nu}}
\otimes V({\varpi}_{\nu})_{(-q_s)^{-\mu}}\,,\\
&& p_{\mu\nu} : V({\varpi}_{\mu})_{(-q_s)^{-\nu}}
\otimes V({\varpi}_{\nu})_{(-q_s)^{\mu}} 
\longrightarrow V({\varpi}_{\mu+\nu}).
\endeqn
\enprop

\subsection{Normalized R-matrices}
Let us calculate $R$-matrices 
between a fundamental representation 
and the vector representation of $U'_q(C_{n}^{(1)})$. 
First recall that we have the following decomposition 
as $U_q(C_{n})$-modules;

\begin{equation}
V(\varpi_{k})\otimes V(\varpi_{1})
=V(\varpi_{k}+\varpi_{1})\oplus V(\varpi_{k+1})\oplus V(\varpi_{k-1}).
\end{equation}
Here $V(\varpi_0)$ is understood to be the trivial representation
and $V(\varpi_{n+1})$ to be $0$.
Therefore the $R$-matrix $\Rn_{k1}(x,y)$ : $V(\varpi_{k})_{x}
\otimes V(\varpi_{1})_{y} \longrightarrow V(\varpi_{1})_{y}\otimes 
V(\varpi_{k})_{x}$ can be written as 
$\Rn_{k1}(x,y)=\Pr_{\varpi _{k+1}+\varpi _{l-1}}
\oplus \gamma_{1}(y/x) \Pr_{\varpi _{k+1}} \oplus \gamma_{2} (y/x) 
\Pr_{\varpi _{k-1}}$, where $\Pr_{\varpi}$ is a $U_q({C_{n}})$-linear 
projection from $V(\varpi_{k})\otimes V(\varpi_{1})$ to 
$V(\varpi)$ in $V(\varpi_{1})\otimes V(\varpi_{k})$ with $\varpi 
=\varpi_{k}+\varpi_{1}$,  $\varpi_{k+1}$ or  $\varpi_{k-1}$.

Let $u_{i}$ and $u'_{i}$ $(i=0,\ 1,\ 2)$ be highest-weight vectors in the
$U_q(C_{n})$-modules $V(\varpi_{k})\otimes V(\varpi_{1})$ and 
$V(\varpi_{1})\otimes V(\varpi_{k})$ with highest weights $\varpi_{k}+
\varpi_{1}\ (i=0)$,  $\varpi_{k+1}\ (i=1)$,  $\varpi_{k-1}\ (i=2)$. 
Remark that if $k=n$ we ignore $\varpi_{k+1},\ u_{1},\ u'_{1}$ and 
$\gamma_{1}(y/x)$. 
We set $Q_{1}=f_{0}f_{1}\cdots f_{n-1}f_{n}f_{n-1}\cdots f_{k+1}$ and 
$Q_{2}=f_{0}f_{1}\cdots f_{k-1}$. 
Then $Q_{i}u_{i}$ is proportional to $u_{0}$ because its weight is 
$\varpi_{k}+\varpi_{1}$. Let us first determine 
$\gamma_{1}$, assuming that $k\neq n$.

The following lemma is by direct calculation and we leave it to the reader.
In the statement $G^{(low)}$ means the lower global base
(cf. \cite{Cry,K}).
\begin{lemma}
Let $b$ be an element of $V(\varpi_{k})\otimes V(\varpi_{1})$ which is 
a tensor product of two lower global bases of $V(\varpi_{k})$ and 
$V(\varpi_{1})$ and has the weight $\varpi_{k+1}$. 
Then $Q_{1}b\neq 0$ if and only if 
$b=b_{1}:=G^{(low)}(1,\ldots,k)\otimes G^{(low)}(k+1)$
 or $b=b_{2}:=G^{(low)}(2,\ldots,k+1)\otimes G^{(low)}(1)$. 
Moreover $Q_{1}b_{1}=q^{-1}y^{-1}u_{0}$ and $Q_{1}b_{2}=x^{-1}u_{0}$, 
where we set $u_{0}=G^{(low)}(1,\ldots,k)\otimes G^{(low)}(1)$.
\end{lemma}

\begin{lemma}
If we write $u_{1}=b_{1} + \sum _{b\neq b_{1}} a_{b}b$, 
where $b$ runs over the set of tensor products of two lower global bases, 
then $a_{b_{2}}=(-q_s)^{k}$.
\end{lemma}

\begin{pf}
There are relations
\eqn
&&e_{i}(G^{(low)}(1,\ldots,\overset{\wedge}{i+1},\ldots,k+1)
\otimes G^{(low)}(i+1) \\
&&\qquad
- q_s G^{(low)}(1,\ldots,\overset{\wedge}{i},\ldots,k+1)\otimes G^{(low)}(i))
=0\quad\mbox{for $1\leq i \leq k$.}
\endeqn
It follows that $a_{b_{2}}=(-q_s)^{k}$.
\end{pf}

By these lemmas we have
$Q_{1}u_{1}=(q_s^{-1}y^{-1}+(-q_s)^{k}x^{-1})u_{0}$ in $V(\varpi_{k})\otimes 
V(\varpi_{1})$.\par

 Similarly we obtain the following two lemmas.
\begin{lemma}
Let $b$ be an element of $V(\varpi_{1})\otimes V(\varpi_{k})$ which is 
a tensor product of two lower global bases of $V(\varpi_{1})$ and 
$V(\varpi_{k})$ and has the weight $\varpi_{k+1}$. 
Then $Q_{1}b\neq 0$ if and only if $b=b'_{1}:=G^{(low)}(1)\otimes 
G^{(low)}(2,\ldots,k+1)$ or $b=b'_{2}:= G^{(low)}(k+1)\otimes 
G^{(low)}(1,\ldots,k)$. 
Moreover $Q_{1}b'_{1}=q_s^{-1}x^{-1}u'_{0}$ and $Q_{1}b'_{2}=y^{-1}u'_{0}$, 
where we set $u'_{0}=G^{(low)}(1)\otimes G^{(low)}(1,\ldots,k)$.
\end{lemma}

\begin{lemma}
If we write $u'_{1}=b'_{1} + \sum _{b\neq b'_{1}} a_{b}b$, where $b$ runs
over the set of tensor products of two lower global bases, then $a_{b'_{2}}=
(-q_s)^{k}$.
\end{lemma}

By these lemmas we have in $V(\varpi_{1})\otimes V(\varpi_{k})$
\begin{equation}  
Q_{1}u'_{1}=(q_s^{-1}x^{-1}+(-q_s)^{k}y^{-1})u'_{0}.
\end{equation}

Therefore we have
\begin{equation}
\gamma_{1}=\frac{x-(-q_s)^{k+1}y}{y-(-q_s)^{k+1}x}.
\end{equation}
Next let us determine $\gamma_{2}$. 
For brevity, we assume that $k\neq 1$ in the following four lemmas,
and $G^{(up)}$ means the lower global base
(cf. \cite{K}).
\begin{lemma}
Let $b$ be an element of $V(\varpi_{k})\otimes V(\varpi_{1})$ which is 
a tensor product of two upper global bases of $V(\varpi_{k})$ and 
$V(\varpi_{1})$ and has the weight $\varpi_{k-1}$. 
Then $Q_{2}b\neq 0$ if and only if 
$b=b_{3}:=G^{(up)}(1,\ldots,k)\otimes 
G^{(up)}(\overset{-}{k})$ or 
$b=b_{4}:=G^{(up)}(2,\ldots,k, \overset{-}{k})\otimes G^{(up)}(1)$.
Moreover $Q_{2}b_{3}=q_s^{-1}y^{-1}u_{0}$ and $Q_{2}b_{4}=q_sx^{-1}[2]_{1}u_{0}$, 
where we set $u_{0}=G^{(up)}(1,\ldots,k)\otimes G^{(up)}(1)$.
\end{lemma}

This lemma is by direct calculation and we leave it to the reader. 

\begin{lemma}
If we write $u_{2}=b_{3} + \sum _{b\neq b_{3}} a_{b}b$, where $b$ runs 
over the set of
tensor products of two upper global bases, then $a_{b_{4}}=
-(-q_s)^{2n-k+1}/[2]_{k-1}$.
\end{lemma}

\begin{pf}
There are relations.\par
\eqn
&&e_{i}(G^{(up)}(1,\ldots,k-1,\ i)\otimes G^{(up)}(\overline{i}) \\
&&\qquad - q_s G^{(up)}(1,\ldots,k-1,\ i+1)\otimes G^{(up)}(\overline{i+1}))=0
\quad\mbox{for $i=k,\ldots,n-1$},
\endeqn
\eqn
&&e_{n}(G^{(up)}(1,\ldots,k-1,\ n)\otimes G^{(up)}(\overline{n}) \\
&&\qquad- q_s^{2} G^{(up)}(1,\ldots,k-1,\ n)\otimes G^{(up)}(\overline{n}))=0,
\endeqn
\eqn
&&e_{i}(G^{(up)}(1,\ldots,k-1,\ \overline{i+1})\otimes G^{(up)}(i+1)  \\
&&\qquad- q_s G^{(up)}(1,\ldots,k-1,\ \overline{i})\otimes G^{(up)}(i))=0
\quad\mbox{for $i=k,\ldots,n-1$,}
\endeqn\eqn
&&e_{k-1}([2]_{k-1}G^{(up)}(1,\ldots,k-1,\ n)\otimes G^{(up)}(\overline{n}) \\ 
&&\qquad- q_s G^{(up)}(1,\ldots,k-1,\ n)\otimes G^{(up)}(\overline{n}))=0,
\endeqn\eqn
&&e_{i}(G^{(up)}(1,\ldots,\overset{\wedge}{i+1},\ldots,k,\ \overline{k})
\otimes G^{(up)}(i+1) \\ 
&&\qquad- q_s G^{(up)}(1,\ldots,\overset{\wedge}{i},\ldots,k,\ \overline{k})
\otimes G^{(up)}(i))=0 \quad\mbox{for $i=1,\ldots,k-2$.}
\endeqn
It follows that $a_{b_{4}}=-(-q_s)^{2n-k+1}/[2]_{k-1}$.
\end{pf}

By these lemmas we have in $V(\varpi_{k})\otimes V(\varpi_{1})$
\begin{equation}
Q_{2}u_{2}=(q_s^{-1}y^{-1}+(-q_s)^{2n-k+2}x^{-1})u_{0}. 
\end{equation}

Similarly we obtain the following two lemmas for $V(\varpi_{1})
\otimes V(\varpi_{k})$.

\begin{lemma}
Let $b$ be an element of $V(\varpi_{1})\otimes V(\varpi_{k})$ which is 
a tensor product of two upper global bases of $V(\varpi_{1})$ and 
$V(\varpi_{k})$ and has the weight $\varpi_{k-1}$. 
Then $Q_{2}b\neq 0$ if and only if $b=b'_{3}:= G^{(up)}(1)\otimes G^{(up)}(2,\ldots,k, \overline{k})$ or $b=b'_{4}:=G^{(up)}(\overline{k})\otimes 
G^{(up)}(1,\ldots,k)$.
Moreover $Q_{2}b'_{3}=q_s^{-1}x^{-1}[2]_{1}u'_{0}$ and 
$Q_{2}b'_{4}=y^{-1}u'_{0}$, where we
set $u'_{0}=G^{(up)}(1)\otimes G^{(up)}(1,\ldots,k)$.
\end{lemma}

\begin{lemma}
If we write $u'_{2}=b'_{3} + \sum _{b\neq b'_{3}} a_{b}b$, where $b$ run
over the set of tensor products of two upper global bases, then $a_{b'_{4}}=
(-q_s)^{2n-k+2}[2]_{k-1}$.
\end{lemma}

By these lemmas we have  in $V(\varpi_{1})\otimes V(\varpi_{k})$
\begin{equation} 
Q_{2}u'_{2}=(q_s^{-1}x^{-1}+(-q_s)^{2n-k+2}y^{-1})(q_s+q_s^{-1})u'_{0}.
\end{equation}

Therefore up to a multiple of an element of $k$ we have
\begin{equation}
\gamma_{2}=\frac{x-(-q_s)^{2n-k+3}y}{y-(-q_s)^{2n-k+3}x}.
\end{equation}

It is easy to check that this expression for $\gamma_{2}$ still holds 
even if $k=1$. So we obtain the following result. 

\begin{thm}
The normalized $R$-matrix is given by
\eqn
\Rn_{k1}(z)=
\left\{
\begin{array}{l}
\Pr_{\varpi _{k}+\varpi _{1}}
+\dfrac{1-(-q_s)^{k+1}z}{z-(-q_s)^{k+1}} \Pr_{\varpi _{k+1}}
+\dfrac{1-(-q_s)^{2n-k+3}z}{z-(-q_s)^{2n-k+3}} \Pr_{\varpi _{k-1}}
\\\noalign{\vspace{5pt}}
\begin{array}{ll}
 &\quad\mbox{if $1\leq k<n$,}\\
\Pr_{\varpi _{k}+\varpi _{1}}
+\dfrac{1-(-q_s)^{n+3}z}{z-(-q_s)^{n+3}}\Pr_{\varpi _{k-1}}
&\quad\mbox{if $k=n$.}
\end{array}\end{array}
\right.
&&
\endeqn
\end{thm}
Hence we have

\eq\label{d:C}
\\\nonumber
&&d_{1,k}(z)=d_{k1}(z)=
\begin{cases}
(z-(-q_s)^{k+1})(z-(-q_s)^{2n+3-k})&\text{if $1\leq k<n$,}\\
z-(-q_s)^{n+3}&\text{if $k=n$.}
\end{cases}
\endeq
We give the explicit form
of $R$-matrix for the vector representation.
\Prop
For $b_1,b_2\in B(\varpi_1)$ we have
\eqn
&&\Rn_{11}(z)(b_1\otimes b_2)=
\begin{cases}
b_1\otimes b_2&\text{if $b_1=b_2$,}\\
\dfrac{(1-q_s^2)z^{\delta(b_2\prec b_1)}}{z-q_s^2}\,b_1\otimes b_2
+\dfrac{q_s(z-1)}{z-q_s^2}\,b_2\otimes b_1
&\text{if $b_1\not= b_2,\ol{b_2}$.}\\
\end{cases}
\endeqn
For $1\le a\le n$ we have
\eqn
&&\Rn_{11}(z)(a\otimes \ol a)=
\dfrac{1-q_s^2}{z-q_s^2}\,a\otimes \ol a
+\sum_{k=1}^n\dfrac
{(-q_s)^{a+k}(1-q_s^2)(z-1)}{(z-q_s^2)(z-(-q_s)^{2n+2})}\,k\otimes \ol k\\
&&\quad\quad-\sum_{k>a}\dfrac
{(-q_s)^{2n+a-k+2}(1-q_s^2)(z-1)}{(z-q_s^2)(z-(-q_s)^{2n+2})}\,\ol k\otimes k
+\dfrac{q_s^{2}(z-1)(z-(-q_s)^{2n})}{(z-q_s^2)(z-(-q_s)^{2n+2})}\,\ol a\otimes a\\
&&\qquad\qquad\qquad-\sum_{k<a}\dfrac
{(-q_s)^{a-k}(1-q_s^2)z(z-1)}{(z-q_s^2)(z-(-q_s)^{2n+2})}\,\ol k\otimes k\,,\\
&&\Rn_{11}(z)(\ol a\otimes a)=
-\sum_{k<a}
\dfrac{(-q_s)^{2n-a+k+2}(1-q_s^2)(z-1)}{(z-q_s^2)(z-(-q_s)^{2n+2})}
\,k\otimes \ol k\\
&&\qquad
+\dfrac{q_s^{2}(z-1)(z-(-q_s)^{2n})}{(z-q_s^2)(z-(-q_s)^{2n+2})}\,a\otimes\ol a
-\sum_{k>a}\dfrac
{(-q_s)^{k-a}(1-q_s^2)z(z-1)}{(z-q_s^2)(z-(-q_s)^{2n+2})}\,k\otimes\ol k\\
&&\qquad+\sum_{k=1}^n\dfrac
{(-q_s)^{2n-a-k+2}(1-q_s^2)z(z-1)}{(z-q_s^2)(z-(-q_s)^{2n+2})}\,\ol k\otimes k
+\dfrac{(1-q_s^2)z}{z-q_s^2}\,\ol a\otimes a\,.\\
\endeqn
\enprop

The general 
$d_{ij}$ with $i,j\not=1$ 
will be calculated at the end of this section with
the aid of the universal $R$-matrices.

\subsection{Universal $R$-matrices}

We shall calculate the universal $R$-matrices.
By (\ref{gen:a}) and (\ref{d:C}), we have
\eq
&&a_{1k}(z)=a_{k1}(z)=
q_s\dfrac{\ip{k-1}\ip{2n+1-k}\ip{2n+3+k}\ip{4n+5-k}}
{\ip{k+1}\ip{2n+3-k}\ip{2n+1+k}\ip{4n+3-k}}\,.
\endeq
Here we employed the notation
\eq\label{ip}
&&\ip{m}=((-q_s)^mz;q_s^{4n+4})_\infty\,.
\endeq
Now we shall calculate $a_{kl}(z)$ for $l\le k$.
Consider the commutative diagram
\eq\label{dia:r}
\\
&&\nonumber
\begin{CD}
V(\varpi_{k})\otimes V(\varpi_{l-1})_{(-q_s)^{-1}z}
\otimes V(\varpi_{1})_{(-q_s)^{l-1}z}
@>{\psi}>>  V(\varpi_{k})\otimes V(\varpi_{l})_{z}\\
\begin{CD}
   @V{f}VV\\
    V(\varpi_{l-1})_{(-q_s)^{-1}z}\otimes V(\varpi_{k})\otimes
     V(\varpi_{1})_{(-q_s)^{l-1}z}\\
    @V{g}VV\\
\end{CD}
@.
@V{h}VV
\\
V(\varpi_{l-1})_{(-q_s)^{-1}z}
\otimes V(\varpi_{1})_{(-q_s)^{l-1}z}
\otimes V(\varpi_{k})
@>>{\psi '}> V(\varpi_{l})_z\otimes V(\varpi_{k})\,.
\end{CD}
\endeq
Here
\eqn
&&\psi=V(\varpi_k)\otimes (p_{l-1,1})_z,
\quad\psi'=(p_{l-1,1})_z\otimes V(\varpi_k),\\
&&f=\Ru_{k,l-1}((-q_s)^{-1}z)\otimes V(\varpi_{1})_{(-q_s)^{l-1}z},\\
&&g=V(\varpi_{l-1})_{(-q_s)^{-1}z}
\otimes \Ru_{k1}((-q_s)^{l-1}z)\ \text{and }
h=\Ru_{kl}(z).
\endeqn
We have
\eqn
p_{l-1,1}(G(1\ldots,l-1)\otimes G(l))&=&G(1,\ldots,l)\,,\\
\Rn_{k,1}(z)(G(1,\ldots,k)\otimes G(l))
&=&G(l)\otimes G(1,\ldots,k).
\endeqn
Chasing the vector
$G(1,\ldots,k)\otimes G(1,\ldots,l-1)\otimes G(l)$
of 
$
V(\varpi_{k})\otimes V(\varpi_{l-1})_{(-q_s)^{-1}z}
\otimes V(\varpi_{1})_{(-q_s)^{l-1}z}$
in the diagram \ref{dia:r},
we obtain the recurrence relation

\eqn
&&
a_{kl}(z)=a_{k,l-1}((-q_s)^{-1}z)a_{k,1}((-q_s)^{l-1}z).
\endeqn
Solving this, and noticing $a_{kl}=a_{lk}$,
we obtain the following result.
\Prop
For $k,l\in I_0=\{1,\ldots,n\}$, we have
\eq
&&a_{kl}(z)=
q_s^{\min(k,l)}\dfrac
{\ip{|k-l|}\ip{2n+2-k-l}\ip{2n+2+k+l}\ip{4n+4-|k-l|}}
{\ip{k+l}\ip{2n+2-k+l}\ip{2n+2+k-l}\ip{4n+4-k-l}}\,.\nn
\endeq
\enprop
Here we used the notation $\ip{m}=((-q_s)^mz;q_s^{4n+4})_\infty$.

\subsection{Denominators of normalized $R$-matrices}
In this subsection we shall prove

\Prop\label{res:C}
For $1\le k,l\le n$, we have
\eq\label{C:d}
\\
&&\nonumber
d_{kl}(z)=
\Prod_{i=1}^{\min(k,l,n-k,n-l)}(z-(-q_s)^{|k-l|+2i})
\Prod_{i=1}^{\min(k,l)}(z-(-q_s)^{2n+2-k-l+2i})\,.
\endeq
\enprop
This is already proved in the case $l=1$.
The case $k=l=n$ is proved in \cite[Proposition 4.2.6]{(KMN)^2_2}.
We shall prove this proposition by reduction to those cases.
Let $D_{kl}(z)$ be the right hand side of (\ref{C:d}).

By (\ref{reverse}), we may assume that $k\ge l$.
First let us show that $d_{kl}(z)$
is a multiple of $D_{kl}(z)$.
In order to see this, 
by using Corollary \ref{Rn:pole},
it is enough to show that
$V(\varpi_k)\otimes V(\varpi_l)_a$ is reducible
for any root $a$ of $D_{kl}(z)$.
For $1\le i\le n-k,l$, we have
\eqn
&&\begin{CD}
V(\varpi_k)\otimes V(\varpi_l)_{(-q_s)^{k-l+2i}}\\
@V{V(\varpi_k)\otimes (i_{i,l-i})_{(-q_s)^{k-l+2i}}}VV\\
V(\varpi_k)\otimes V(\varpi_i)_{(-q_s)^{k+i}}\otimes
V(\varpi_{l-i})_{(-q_s)^{k-l+i}}\\
@V{(p_{ki})_{(-q_s)^i}\otimes V(\varpi_{l-i})_{(-q_s)^{k-l+i}}}VV\\
V(\varpi_k)_{(-q_s)^i}\otimes V(\varpi_{l-i})_{(-q_s)^{k-l+i}}
\end{CD}
\endeqn
Here $V(\varpi_0)$ is understood to be the trivial representation.
Then one can easily see that
the composition is not zero but $u_k\otimes u_l$
is sent to zero.
Hence $V(\varpi_k)\otimes V(\varpi_l)_{(-q_s)^{k-l+2i}}$
is reducible.
Similarly for $1\le i\le l$, let us consider
\eqn
&&\begin{CD}
V(\varpi_k)\otimes V(\varpi_l)_{(-q_s)^{2n+2-k-l+2i}}\\
@V{i_{k-i,i}\otimes (i_{i,l-i})_{(-q_s)^{2n+2-k-l+2i}}}VV\\
V(\varpi_{k-i})_{(-q_s)^i}\otimes V(\varpi_i)_{(-q_s)^{i-k}}\otimes
V(\varpi_i)_{(-q_s)^{2n+2-k+i}}\otimes
V(\varpi_{l-i})_{(-q_s)^{2n+2-k-l+i}}\\
@V{V(\varpi_{k-i})_{(-q_s)^i}\otimes \tr\otimes
V(\varpi_{l-i})_{(-q_s)^{2n+2-k-l+i}}}VV\\
V(\varpi_{k-i})_{(-q_s)^i}\otimes V(\varpi_{l-i})_{(-q_s)^{2n+2-k-l+i}}
\end{CD}
\endeqn
In this case also, the composition is not zero but $u_k\otimes u_l$
is sent to zero.
Hence $V(\varpi_k)\otimes V(\varpi_l)_{(-q_s)^{2n+2-k-l+2i}}$ is reducible.

By (\ref{eq:univ}), we have
$$a_{kl}(z)\ake a_{kl}((-q_s)^{-(2n+2)}z)\equiv
\dfrac{d_{kl}(z)}{d_{kl}((-q_s)^{2n+2}z^{-1})}\quad\mod k[z,z^{-1}]^\times\,.$$
Hence we obtain
\eq
d_{kl}(z)=D_{kl}(z)\psi_{kl}(z)
\endeq
for a polynomial $\psi_{kl}(z)$ satisfying
\eq\label{sym:psi}
\psi_{kl}(z)\equiv\psi_{kl}((-q_s)^{2n+2}z^{-1})
\quad\mod k[z,z^{-1}]^\times.
\endeq

Now we shall use the following lemma.
\Lemma
Let $V'$, $V''$, $V$ and $W$ be irreducible $\Us$-modules.
Assume that there is a surjective morphism
$V'\otimes V''\to V$.
Then
$$
\dfrac{d_{W,V'}(z)d_{W,V''}(z)a_{W,V}(z)}
{d_{W,V}(z)a_{W,V'}(z)a_{W,V''}(z)}
\quad\text{and}\quad
\dfrac{d_{V',W}(z)d_{V'',W}(z)a_{V,W}(z)}
{d_{V,W}(z)a_{V',W}(z)a_{V'',W}(z)}
$$
are in $k[z,z^{-1}]$.
\enlemma
\proof
In a commutative diagram
\eqn
\begin{CD}
W\otimes V'_z\otimes V''_z@>>>W\otimes V_z\\
\begin{CD}
@V{R'(z)\otimes V''_z}VV\\
V'_z\otimes W\otimes V''_z\\
@V{V'_z\otimes R''(z)}VV
\end{CD}
@.
\begin{CD}
@V{R(z)}VV
\end{CD}\\
V'_z\otimes V''_z\otimes W
@>>>V_z\otimes W\,,
\end{CD}
\endeqn
if $R'(z)$ and $R''(z)$ do not have poles, then so is $R(z)$.
To see the first assertion, it is enough to apply this to
$R'(z)=d_{W,V'}(z)\Rn_{W,V'}(z)$,
$R''(z)=d_{W,V''}(z)\Rn_{W,V''}(z)$
and 
$$R(z)=\dfrac{d_{W,V'}(z)d_{W,V''}(z)a_{W,V}(z)}
{d_{W,V}(z)a_{W,V'}(z)a_{W,V''}(z)}\Rn_{W,V}(z)\,.$$
The second assertion can be proved similarly.
\QED

We shall prove $\psi_{kl}(z)\equiv 1$ $\mod k[z,z^{-1}]^\times$.

\medskip
\noindent
{\bf Case $k+l\le n$.}
\quad
We prove this by the induction on $l$.
If $l=1$, it is already proved.
If $l>1$ then applying the lemma above to
$V(\varpi_{l-1})_{(-q_s)^{-1}}\otimes 
V(\varpi_1)_{(-q_s)^{l-1}}\to V(\varpi_l)$,
we have
\eqn
&&
\dfrac{d_{k,l-1}((-q_s)^{-1}z)\ake d_{k,1}((-q_s)^{l-1}z)\ake a_{k,l}(z)}
{d_{k,l}(z)\ake a_{k,l-1}((-q_s)^{-1}z)\ake a_{k,1}((-q_s)^{l-1}z)}
\in k[z,z^{-1}].
\endeqn
Since 
\eqn
&&
\dfrac{D_{k,l-1}((-q_s)^{-1}z)\ake d_{k,1}((-q_s)^{l-1}z)\ake a_{k,l}(z)}
{D_{k,l}(z)\ake a_{k,l-1}((-q_s)^{-1}z)\ake a_{k,1}((-q_s)^{l-1}z)}\equiv 1,
\endeqn
$\psi_{k,l-1}(z)\equiv 1$ implies $\psi_{kl}(z)\equiv 1$.

\medskip
\noindent
{\bf Case $k+l>n$}\quad
We shall first reduce the assertion to the $k=n$ case.
For $k<n$ consider a surjection
\eqn
&&V(\varpi_{k+1})_{(-q_s)^{-1}}\otimes V(\varpi_1)_{(-q_s)^{2n+1-k}}
\to V(\varpi_k)
\endeqn
given by the composition
\eqn
&&V(\varpi_{k+1})_{(-q_s)^{-1}}\otimes V(\varpi_1)_{(-q_s)^{2n+1-k}}
\to
V(\varpi_{k})\otimes V(\varpi_1)_{(-q_s)^{-1-k}}
\otimes V(\varpi_1)_{(-q_s)^{2n+1-k}}\\
&&\phantom{
V(\varpi_{k+1})_{(-q_s)^{-1}}\otimes V(\varpi_1)_{(-q_s)^{2n+1-k}}
\to
V(\varpi_{k})\otimes V(\varpi_1)_{(-q_s)^{-1-k}}
\otimes }
\to V(\varpi_{k})\,.
\endeqn
We have
\eqn
&&\dfrac{D_{k+1,l}((-q_s)z)\ake d_{1,l}((-q_s)^{k-2n-1}z)\ake a_{kl}(z)}
{D_{kl}(z)\ake a_{k+1,l}((-q_s)z)\ake a_{1,l}((-q_s)^{k-2n-1}z)}
\equiv z-(-q_s)^{4n+4-k-l}.
\endeqn
Hence $\psi_{k+1,l}(z)\equiv1$ implies that $\psi_{kl}(z)$ is a divisor
of $z-(-q_s)^{4n+4-k-l}$. Then (\ref{sym:psi}) implies that 
$\psi_{kl}(z)\equiv 1$.
Hence, the descending induction on $k$
reduces the problem to the $k=n$ case.
\hb
We have
\eqn
&&
\dfrac{D_{k,l-1}((-q_s)^{-1}z)\ake d_{k,1}((-q_s)^{l-1}z)\ake a_{k,l}(z)}
{D_{k,l}(z)\ake a_{k,l-1}((-q_s)^{-1}z)\ake a_{k,1}((-q_s)^{l-1}z)}
\equiv z-(-q_s)^{2n+2-k+l}.
\endeqn
Hence by the similar argument to $k+l\le n$ case,
$\psi_{k,l-1}(z)\equiv 1$ implies that $\psi_{kl}(z)$ is a divisor
of $z-(-q_s)^{2n+2-k+l}$. Hence if $l\not=k=n$ then
we can reduce the $l$ case to the $l-1$ case.
This completes the proof of Proposition \ref{res:C}.

\bibliographystyle{unsrt}

\end{document}